\documentstyle{amltd}
\begin{document}
\annalsline{157}{2003}
\received{May 10, 2001}
\startingpage{847}
\def\bye{\end{document}}
 \font\tenrm=cmr10
\def\ritem#1{\item[{\rm #1}]}
 \input amssym.def
\input amssym.tex

\input boxedeps.tex 
\SetepsfEPSFSpecial 
\HideDisplacementBoxes
\def\figin#1#2{
$$
 {\BoxedEPSF{#1.eps scaled
#2}%
}%
$$
\noindent}

\catcode`\@=11
\font\twelvemsb=msbm10 scaled 1100
\font\tenmsb=msbm10
\font\ninemsb=msbm10 scaled 800
\newfam\msbfam
\textfont\msbfam=\twelvemsb  \scriptfont\msbfam=\ninemsb
  \scriptscriptfont\msbfam=\ninemsb
\def\msb@{\hexnumber@\msbfam}
\def\Bbb{\relax\ifmmode\let\next\Bbb@\else
 \def\next{\errmessage{Use \string\Bbb\space only in math
mode}}\fi\next}
\def\Bbb@#1{{\Bbb@@{#1}}}
\def\Bbb@@#1{\fam\msbfam#1}
\catcode`\@=12

 \catcode`\@=11
\font\twelveeuf=eufm10 scaled 1100
\font\teneuf=eufm10
\font\nineeuf=eufm7 scaled 1100
\newfam\euffam
\textfont\euffam=\twelveeuf  \scriptfont\euffam=\teneuf
  \scriptscriptfont\euffam=\nineeuf
\def\euf@{\hexnumber@\euffam}
\def\frak{\relax\ifmmode\let\next\frak@\else
 \def\next{\errmessage{Use \string\frak\space only in math
mode}}\fi\next}
\def\frak@#1{{\frak@@{#1}}}
\def\frak@@#1{\fam\euffam#1}
\catcode`\@=12

\newcommand{\thmref}[1]{Theorem~\ref{#1}}
\newcommand{\secref}[1]{\S\ref{#1}}
\newcommand{\lemref}[1]{Lemma~\ref{#1}}
\title{Rogers-Ramanujan and the\\ Baker-Gammel-Wills (Pad\'{e}) conjecture} 
\shorttitle{The Baker-Gammel-Wills conjecture} 
\acknowledgements{}
\author{D.\ S.\  Lubinsky}
  \institutions{The John Knopfmacher Centre, University of Witwatersrand, Wits, South
Africa and  Georgia Institute of Technology,
Atlanta,  GA.\\
{\eightpoint {\it E-mail address\/}: lubinsky@math.gatech.edu
}}


  \centerline{\ninepoint \it Dedicated to the memory of Israel{\rm ,} Zivia and Ranan Lubinsky}
\vglue12pt
\centerline{\bf Abstract}
\vglue12pt

In 1961, Baker, Gammel and Wills conjectured that for functions $f$
meromorphic in the unit ball, a subsequence of its diagonal Pad\'{e}
approximants converges uniformly in compact subsets of the ball omitting
poles of $f$. There is also apparently a cruder version of the conjecture
due to Pad\'{e} himself, going back to the early twentieth century. We
show here that for carefully chosen $q$ on the unit circle, the
Rogers-Ramanujan continued fraction 
$$
1+\frac{qz|}{|1}+\frac{q^{2}z|}{|1}+\frac{q^{3}z|}{|1}+\cdots 
$$
provides a counterexample to the conjecture. We also highlight some other
interesting phenomena displayed by this fraction.

\section{Introduction}

Let 
$$
f\left( z\right) =\sum_{j=0}^{\infty }a_{j}z^{j}
$$
be a formal power series, with complex coefficients. Given integers $m,n\geq
0$, the $\left( m,n\right) $ Pad\'{e} approximant to $f$ is a rational
function 
$$
\left[ m/n\right] =P/Q
$$
where $P,Q$ are polynomials of degree at most $m,n$ respectively, such that $%
Q$ is not identically $0$, and such that 
\begin{equation}
\left( fQ-P\right) \left( z\right) =O\left( z^{m+n+1}\right) .
\end{equation}
By this last relation, we mean that the coefficients of $%
1,z,z^{2},\ldots,z^{m+n}$ in the formal power series on the left-hand side
vanish. The basic idea is that $\left[ m/n\right] $ is a rational function
with given upper bounds on its numerator and denominator degrees, chosen in
such a way that its Maclaurin series reproduces as many terms as possible in
the power series $f$.

It is easy to see that $\left[ m/n\right] $ exists: we can reformulate (1.1)
as a system of $m+n+1$ homogeneous linear equations in the $\left(
m+1\right) +\left( n+1\right) $ coefficients of the polynomials $P$ and $Q$.
As there are more unknowns than equations, there is a nontrivial solution,
and it is easily seen from (1.1) that $Q$ cannot be identically $0$ in any
nontrivial solution. While $P$ and $Q$ are not separately unique, the ratio 
$\left[ m/n\right] $ is, and this is again an easy consequence of (1.1).

It was C. Hermite, who gave his student Henri Eugene Pad\'{e} the
approximant to study in the 1890's. Although the approximant was known
earlier, by amongst others, Jacobi and Frobenius, it was perhaps Pad\'{e}'s
thorough investigation of the structure of the Pad\'{e} table, namely the
array 
$$
\begin{array}{lllll}
\left[ 0/0\right] & \left[ 0/1\right] & \left[ 0/2\right] & \left[ 0/3\right]
& \ldots \\[4pt]
\left[ 1/0\right] & \left[ 1/1\right] & \left[ 1/2\right] & \left[ 1/3\right]
& \ldots \\ [4pt]
\left[ 2/0\right] & \left[ 2/1\right] & \left[ 2/2\right] & \left[ 2/3\right]
& \ldots \\ [4pt]
\left[ 3/0\right] & \left[ 3/1\right] & \left[ 3/2\right] & \left[ 3/3\right]
& \ldots \\ [4pt]
\vdots & \vdots & \vdots & \vdots & \ddots
\end{array}
$$
that has ensured the approximant being named after him.

Pad\'{e} approximants have been applied in proofs of irrationality and
transcendence in number theory, in practical computation of special
functions, and in analysis of difference schemes for numerical solution of
partial differential equations. However, the application which really
brought them to prominence in the 1960's and 1970's, was in location of
singularities of functions: in various physical problems, for example
inverse scattering theory, one would have a means for computing the
coefficients of a power series $f$. One could use these coefficients to
compute, for example, the $[3/3]$ Pad\'{e} approximant to $f$, and use the
poles of the approximant as predictors of the location of poles or other
singularities of $f$. Moreover, under certain conditions on $f$, which were
often satisfied in physical examples, this process could be theoretically
justified.

In addition to their wide variety of applications, they are also closely
associated with continued fraction expansions, orthogonal polynomials,
moment problems, the theory of quadrature, amongst others. See \cite
{Baker1975} and \cite{BakerGravesMorris1996} for a detailed development of
the theory, and \cite{Brezinski1991} for their history.

One of the fascinating features of Pad\'{e} approximants is the complexity
of their convergence theory. There are power series $f$ with zero radius of
convergence, for which $\left[ n/n\right] \left( z\right) $ converges as $%
n\rightarrow \infty $ to a function single valued and analytic in the
cut-plane ${\Bbb C}\backslash [0,\infty )$. On the other hand, there are
entire functions $f$ for which 
$$
\limsup_{n\rightarrow \infty }\left| [n/n]\left( z\right) \right| =\infty
$$
for all $z\in {\Bbb C}\backslash \left\{ 0\right\} $.

Probably the most important general theorem that applies to functions
meromorphic in the plane is that of Nuttall-Pommerenke. It asserts that if $%
f $ is meromorphic throughout ${\Bbb C}$, and analytic at $0$, then $\left\{
\left[ n/n\right] \right\} _{n=1}^{\infty }$ converges in planar measure.
More generally, this holds if $f$ has singularities of (logarithmic)
capacity $0$, and planar measure may be replaced by capacity. There are much
deeper analogues of this theorem for functions with branchpoints, due to H.
Stahl. Uniform convergence of sequences of Pad\'{e} approximants has been
established for P\'olya frequency series, series of
Stieltjes/Markov/Hamburger, and other special classes. For surveys and
various perspectives on the convergence theory, see \cite{Baker1975}, \cite
{BakerGravesMorris1996}, \cite{Goncar1982}, \cite{Lubinsky1999-1}, \cite
{Lubinsky2001-2}, \cite{Stahl1989}, \cite{Stahl1997-1}, \cite{Stahl1997-2}, 
\cite{Wallin1979}.

Long before the Nuttall-Pommerenke theorem was established, George Baker and
his collaborators observed the phenomenon of spurious poles: several of the
approximants could have poles which in no way were related to those of the
underlying function. However, those poles affected convergence only in a
small neighbourhood, and there were usually very few of these ``bad''
approximants. Thus, one might compute $\left[ n/n\right] ,\;n=1,2,3,\ldots50$,
and find a definite convergence trend in 45 of the approximants, with five of
the 50 approximants displaying pathological behaviour. The curious thing
(contrary to expectation) is that the five bad approximants could be
distributed anywhere in the 50, and need not be the first few. Nevertheless,
after omitting the ``bad'' approximants, one obtained a clear convergence
trend. This seemed to be a characteristic of the Pad\'{e} method, and Baker et al.\
formulated a now famous conjecture \cite{BakerGammelWills1961}. There are
now many forms of the conjecture; we shall concentrate on the following form:%

\nonumproclaim{Baker-Gammel-Wills Conjecture (1961)} 
 Let $f$ be meromorphic in the unit ball{\rm ,} and analytic at 
$0$.  There is an infinite subsequence $\left\{ [n/n]\right\} _{n\in 
{\cal S}}$ of the diagonal sequence  $\left\{ \left[ n/n\right]
\right\} _{n=1}^{\infty }$ that converges uniformly in all compact
subsets of the unit ball omitting poles of $f$.
\endproclaim

Thus, there is an infinite sequence of ``good'' approximants. In the first
form of the conjecture, $f$ was required to have a nonpolar singularity on
the unit circle, but this was subsequently relaxed (cf. \cite[p.~188 ff.]
{Baker1975}). There is also apparently a cruder form of the conjecture due
to Pad\'{e} himself, dating back to the 1900's; the author must thank J.
Gilewicz for this historical information.

The main result of this paper is that the above form of the conjecture is
false, and that a counterexample is provided by a continued fraction of
Rogers-Ramanujan. For $q$ not a root of unity, let 
\begin{equation}
G_{q}\left( z\right) :=\sum_{j=0}^{\infty }\frac{q^{j^{2}}}{\left(
1-q\right) \left( 1-q^{2}\right)\cdots \left( 1-q^{j}\right) }z^{j}
\end{equation}
denote the Rogers-Ramanujan function. Of course, it is at this stage merely
a formal power series. Moreover, let 
\begin{equation}
H_{q}\left( z\right) :=G_{q}\left( z\right) /G_{q}\left( qz\right) .
\end{equation}
When $H_{q}$ has an analytic (or meromorphic) continuation to a region
beyond the domain of definition of $G_{q}$, we denote that continuation by $%
H_{q}$ also. There is the well-known functional relation, which we shall
establish in Section 3: 
\begin{equation}
H_{q}\left( z\right) =1+\frac{qz}{H_{q}\left( qz\right) }.
\end{equation}
Iterating this leads to 
\begin{equation}
H_{q}\left( z\right) =1+\frac{qz}{1+\frac{q^{2}z}{1+\ddots \frac{q^{n}z}{%
H_{q}\left( q^{n}z\right) }}}
\end{equation}
and hence to the formal infinite continued fraction 
\begin{equation}
H_{q}\left( z\right) =1+\frac{qz|}{|1}+\frac{q^{2}z|}{|1}+\frac{q^{3}z|}{|1}%
+\cdots \;.
\end{equation}
(The continued fraction notation used should be self
explanatory.)
For\break $\left| q\right|  < 1$, the continued fraction was considered
independently by L.\ J.\  Rogers and  S.\ Ramanujan in the early part of the
twentieth century.

The truncations of a continued fraction are called its \textit{convergents}.
We shall use the notation 
\begin{equation}
\frac{\mu _{n}}{\nu _{n}}\left( z\right) =1+\frac{qz|}{|1}+\frac{q^{2}z|}{|1}%
+\cdots +\frac{q^{n}z|}{|1},\qquad n\geq 1
\end{equation}
for the $n^{\rm th}$ convergent, to emphasize that it is a rational function with
numerator polynomial $\mu _{n}$ and denominator polynomial $\nu _{n}$. We also set 
$$
\mu _{0}/\nu _{0}:=1.
$$
The continued fraction is said to \textit{converge} if 
$$
\lim_{n\rightarrow \infty }\mu _{n}\left( z\right) /\nu _{n}\left( z\right)
$$
exists.

At least when $G_{q}$ has a positive radius of convergence, it does not
really matter whether we define $H_{q}$ by (1.3) or (1.6), for both have the
same Maclaurin series, so both analytically continue that Maclaurin series
inside their domain of convergence. When $G_{q}$ has zero radius of
convergence, we shall define $H_{q}$ by (1.6).

We shall make substantial use of the fact that the sequence $\left\{ \mu
_{n}/\nu _{n}\right\} _{n=1}^{\infty }$ of convergents includes both the
diagonal sequence $\left\{ \left[ n/n\right] \right\} _{n=1}^{\infty }$ and
the sub-diagonal sequence $\left\{ \left[ n+1/n\right] \right\}
_{n=1}^{\infty }$ to $H_{q}$. So as $n$ increases, the convergents trace a
stair step in the Pad\'{e} table. For a proof of this, see \cite
{BakerGravesMorris1996} or \cite{LorentzenWaadeland1992}.

Our counterexample is contained in: 

\nonumproclaim{Theorem 1.1} Let  
\begin{equation}
q:=\exp \left( 2\pi i\tau \right)
\end{equation}
 where 
\begin{equation}
\tau :=\frac{2}{99+\sqrt{5}}.
\end{equation}
 Then  $H_{q}$  is meromorphic in the unit ball and analytic
at $0$. There does not exist any subsequence of  $\left\{ \mu
_{n}/\nu _{n}\right\} _{n=1}^{\infty }$  that converges uniformly in
all compact subsets of  
$$
{\cal A}:=\left\{ z:\left| z\right| <0.46\right\}
$$
 omitting poles of $H_{q}$. In particular no subsequence
of  $$ \left\{ \left[ n/n\right] \right\} _{n=1}^{\infty }\qquad \hbox{ or } \qquad 
\left\{ \left[ n+1/n\right] \right\} _{n=1}^{\infty }$$  can converge
uniformly in all compact subsets of ${\cal A}$  omitting poles
of  $H_{q}$.
\endproclaim

The crux of the counterexample is that, given any subsequence $\left\{ \mu
_{n}/\nu _{n}\right\} _{n\in {\cal S}}$ of the convergents, there is a
compact subset of ${\cal A}$ not containing any poles of $H_{q}$, such
that infinitely many of the convergents have a pole in the interior of the
compact set. Moreover, there is a limit point of poles in the interior of
that compact set, and uniform convergence is not possible.

There are several limits to our example. We are certain that with sufficient
effort, one may replace $0.46$ above by $\frac{1}{4}+\varepsilon $, for an
arbitrarily small $\varepsilon >0$ and a corresponding $q$ on the unit
circle. However, we cannot go below $\frac{1}{4}$. Indeed, an old theorem of
Worpitzky guarantees that the full sequence of convergents $\left\{ \mu
_{n}/\nu _{n}\right\} _{n=1}^{\infty }$ converges uniformly in compact
subsets of $\left\{ z:\left| z\right| <\frac{1}{4}\right\} $. Thus one can
still look for an example in which no subsequence of the convergents
converges uniformly, or even pointwise, in any neighbourhood of $0$.

Moreover, given any point in the unit ball at which $H_{q}$ is analytic,
there is a neighbourhood of it and a subsequence of the convergents that
converges uniformly in that neighbourhood. So one can also look for an
example without this property. We shall discuss this further in Section 8.

We shall see that for a.e. $q$ on the unit circle (and in particular for the 
$q$ above), $H_{q}$ is meromorphic in the unit ball, with a natural boundary
on the unit circle. Moreover, for a.e. $q$, $G_{q}$ is analytic in the unit
ball, with a natural boundary on the unit circle.

However, given $0<s<\frac{1}{4}$, then for some exceptional $q$, there is
the very striking feature, that $G_{q}$ is analytic in $\left| z\right| <s$,
with a natural boundary on the circle $\left\{ z:\left| z\right| =s\right\} $%
, yet $H_{q}$ defined by (1.3) admits an analytic continuation to at least
the ball centre $0$, radius $\frac{1}{4}$. So somehow, in the division in
(1.3), the natural boundary of $G_{q}$ is cancelled out, as if it were a
removable singularity.

There are other striking features for a.e. $q$: if on a circle centre $0$, $%
H_{q}$ has poles of total multiplicity $\ell $, then in any neighbourhood of
that circle, all convergents $\mu _{n}/\nu _{n}$ with $n$ large enough, have
at least $2\ell $ poles, namely double as many as $H_{q}$.

This paper is organised as follows: in Section 2, we shall state in greater
detail, our results on $G_{q}$, $H_{q}$ and the convergence or divergence
properties of the continued fraction. In Section 3, we shall present some
identities involving the approximants and their proofs. In Section 4, we
shall prove our results on the continued fraction when $q$ is a root of
unity. In Sections 5 and 6, we shall prove the results of Section 2. In
Section 7, we shall prove Theorem 1.1. Finally in Section 8, we shall
discuss some of the implications of this paper.

\section{The continued fraction for $H_{q}$}

We emphasise that the Rogers-Ramanujan c.f.\ (continued fraction) is not the
first candidate we have examined as a possible counterexample to the
Baker-Gammel-Wills conjecture. In the search for a counterexample, basic
hypergeometric, or $q$ series, have been most useful, just as they have
had applications in so many branches of mathematics. What is somewhat
exotic, however, is the range of the parameter $q$. In most studies of $q$%
-series, $\left| q\right| <1$, and sometimes $\left| q\right| >1$. However,
many of the identities persist for $\left| q\right| =1$, and it is in this
range of $q$, that several interesting phenomena and counterexamples in the
convergence theory of Pad\'{e} approximation have been discovered. In other
contexts, the case $\left| q\right| =1$ has also proved to be interesting 
\cite{SpiridonovZhedanov2000}.

In \cite{LubinskySaff1987}, E.\ B.\  Saff and the author investigated the
Pad\'{e} table and continued fraction for the partial theta function 
$$
\sum_{j=0}^{\infty }q^{j(j-1)/2}z^{j}=1+\frac{z|  }{|  1}-\frac{qz|  }{%
|  1}+\frac{q(1-q)z|  }{|  1}-\frac{q^{3}z|  }{|  1}+\frac{%
q^{2}(1-q)z|  }{|  1}\cdots
$$
when $|q|=1$. Subsequently K.\ A.\  Driver and the author \cite
{Driver1991}, \cite{DriverLubinsky1990}, \cite{DriverLubinsky1991}, \cite
{DriverLubinsky1993} undertook a detailed study of the Pad\'{e} table and
continued fraction for the more general Wynn's series \cite{Wynn1967} 
\begin{eqnarray*}
&&\sum_{j=0}^{\infty }\left[\prod_{l=0}^{j-1}(A-q^{l+\alpha })\right]z^{j}; \\
&& \\
&&\sum_{j=0}^{\infty }\frac{z^{j}}{\prod_{l=0}^{j-1}(C-q^{l+\alpha })}; \\
&& \\
&&\sum_{j=0}^{\infty }\left[\prod_{l=0}^{j-1}\frac{A-q^{l+\alpha }}{C-q^{l+\gamma
}}\right]z^{j}.
\end{eqnarray*}
Here $A,C,\alpha $ and $\gamma $ are suitably restricted parameters.

Amongst the interesting features is that some subsequence of the convergents
converges uniformly inside the region of analyticity, so that
Baker-Gammel-Wills is true for these series, while ``most'' subsequences
have poles that cycle around the region of analyticity.

There are at least three aspects of the Rogers-Ramanujan c.f.\ that
distinguish it from Wynn's series in the case where $|q|=1$. Firstly
the functional relation for $H_{q}$, namely 
$$
H_{q}(z)=1+\frac{qz}{H_{q}(qz)}
$$
generates its c.f.\ by repeated application. For Wynn's series, there is not
such a simple relationship between the c.f.\ and the functional equation.
Secondly all the coefficients in the Rogers-Ramanujan c.f.\ have modulus $1$,
whereas a subsequence of the coefficients in the c.f.\ for Wynn's series
converges to $0$. Moreover the latter subsequence is associated with a
subsequence of the convergents to the c.f.\ that converges throughout the
region of analyticity. This already suggests that there may not be a
uniformly convergent subsequence of the convergents for the Rogers-Ramanujan
c.f. Thirdly, in the case where $q$ is a root of unity, all of the Wynn's
series reduce to rational functions, while the Rogers-Ramanujan c.f.\
corresponds to a function with branchpoints.

It is an immediate consequence of Worpitzky's theorem that the c.f.\ (1.6)
converges for $|z|<\frac{1}{4}$, for each $|q|=1$. In fact,
we shall show using standard methods that (1.6) converges for $|z|<(2+|1+q|)^{-1}$. However beyond that circle, standard methods give
very little, because of the oscillatory nature of the continued fraction
coefficients $\left\{ q^{n}\right\} _{n=1}^{\infty }$.

One must obviously distinguish the case where $q$ is a root of unity, as the
power series coefficients of $G_{q}$ are not even defined in this case.  Then, rather than defining $H_{q}$ by (1.3), we
shall define it as the function corresponding to the continued fraction (1.6). Using standard
results for periodic c.f.'s, we shall prove in Section 4, the following:%

\nonumproclaim{Theorem 2.1} Let $\ell \geq 1$  and  $q$  be a primitive  $\ell^{\rm th}$ root of unity. Let 
\begin{equation}
{\cal L}:=\left\{z\in {\Bbb C}:z^{\ell }\in \left(-\infty ,-\frac{1}{4}\right]\right\}.
\end{equation}
 There exists a set  ${\cal P}$  of at most  $(\ell
-1)(2\ell -1)/2$  points such that for  $z\in {\Bbb C}\backslash (%
{\cal L}\cup {\cal P})${\rm ,}
\begin{equation}
\lim_{n\rightarrow \infty }\frac{\mu _{n}(z)}{\nu _{n}(z)}=\frac{\mu _{\ell
-1}(z)-\frac{1}{2}+\sqrt{z^{\ell }+\frac{1}{4}}}{\nu _{\ell -1}(z)}=:H_{q}(z).
\end{equation}
 Here the branch of  $\sqrt{}$ is the principal one{\rm ,}
analytic in ${\Bbb C}\backslash (-\infty ,0]$  and positive in  $%
\left( 0,\infty \right) $.
\endproclaim

Of course, ${\cal L}$ consists of $\ell $ distinct rays with an angle of $%
2\pi /\ell $ between successive rays, extending from the $\ell $ values of $%
(-\frac{1}{4})^{1/\ell }$ out to $\infty $. So the c.f.\   chooses the most
natural choice for the branchcuts; see \cite{Stahl1989}, \cite{Stahl1997-1}
for the ways that continued fractions and Pad\'{e} approximants choose
branchcuts in far more general situations.

The set ${\cal P}$ contains the poles of $H_{q}$, that is, the at most $%
(\ell -1)/2$ zeros of $\nu _{\ell -1}$, which need not lie on the branchcuts
contained in the set ${\cal L}$. For example, if $\ell =5$, $\nu _{\ell
-1}(z)=(qz-1)(z-1)$ has zeros at $1$ and $\overline{q}$, which are not in $%
{\cal L}$. Also, ${\cal P}$  contains additional points that arise in
applying the standard theorems on periodic continued fractions. We have not
been able to determine if these additional points are really points of
divergence, or to determine where they lie. In all probability, our bound of 
$(\ell -1)(2\ell -1)/2$ on the number of points in ${\cal P}$ is too
large.

Next, we turn to the more difficult case where $q$ is not a root of unity.
Clearly the series $G_{q}$ of (1.2) at least has well-defined coefficients
if $q$ is not a root of unity, and its radius of convergence is 
\begin{equation}
R(q):=\liminf_{j\rightarrow \infty }\left|\prod_{k=0}^{j-1}(1-q^{k})\right|^{1/j}.
\end{equation}
It was essentially proved in \cite{HardyLittlewood1946} (and we shall
reproduce the proof in Lemma 6.2) that 
\begin{equation}
R(q)=\liminf_{j\rightarrow \infty }\left|1-q^{j}\right|^{1/j}.
\end{equation}
If we write $q=e^{2\pi i\tau }$, this is readily reformulated in terms of
the diophantine approximation properties of $\tau $. Since $|1-q^{j}|=2|\sin [\pi (j\tau -k)]|$ for any integer $k$, we see that 
\begin{equation}
R(q)=\liminf_{j\rightarrow \infty }\Vert j\tau \Vert^{1/j},
\end{equation}
where $\Vert x\Vert $ denotes the distance from $x$ to the nearest
integer.

It is known that $R(q)=1$ for ``most'' $q$. Indeed the set 
\begin{equation}
{\cal G}:=\{q:R(q)<1\}
\end{equation}
has linear measure $0$, Hausdorff dimension $0$, and even logarithmic
dimension $2$ \cite{Lubinsky1985}. G.\ Petruska has shown \cite{Petruska1992}
that the related quantity 
$$
\limsup_{j\rightarrow \infty }\left|\prod_{k=0}^{j-1}(A-q^{k})\right|^{1/j} 
$$
may assume any value in $[0,1]$ as $A$ and $q$ range over the unit circle.
Using his results, we can easily show that $R(q)$ may assume any value in $%
[0,1]$. Curiously enough, the radius of convergence $R(q)$ of $G_{q}$ need
not coincide with the radius of meromorphy of $H_{q}$, that is, the largest
circle centre $0$ inside which $H_{q}$ may be meromorphically continued. On
the boundary of that circle, we show that $H_{q}$ has a natural boundary:%

\nonumproclaim{Theorem 2.2} 
 Let $|q|=1${\rm ,} and assume that  $q$ is not a
root of unity. Let $\rho (q)$ denote the radius of meromorphy of  $%
H_{q}$. Then 
\vglue4pt
{\rm (a) }$H_{q}$ has a natural boundary on the circle  $\{z:|z|=\rho (q)\}$ and  
\begin{equation}
1\geq \rho (q)\geq \max \left\{R(q),\frac{1}{2+|1+q|}\right\}\geq \frac{1}{4}.
\end{equation}
\vglue4pt
{\rm (b) }$G_{q}$ has a natural boundary on the circle $\{z:|z|=R(q)\}.$ Moreover{\rm ,} as $q$  ranges over the unit circle{\rm,} $%
R\left( q\right) $ may assume any value in  $[0,1]$.

\vglue4pt
{\rm (c) } For $q\notin {\cal G},R(q)=\rho (q)=1$. In particular{\rm ,} this
is true for {\rm a.e.}  $q$.
\endproclaim

We are not sure if $\rho (q)$ may assume values $<1$, but are inclined to
believe that always $\rho (q)=1$. At least for ``most'' $q$, the above
result asserts that $H_{q}$ is given by (1.3) inside its radius of
meromorphy.

We are also interested in how $H_{q}$ varies as $q$ does, especially near
roots of unity, as the branchcuts of $H_{q}$ should then attract poles and
zeros of the ``nearby'' meromorphic $H_{q}$. The following result partly
justifies the latter: 

\nonumproclaim{Theorem 2.3} Let $|q_{k}|=1,k\geq 1${\rm ,} and assume that   
\begin{equation}
\lim_{k\rightarrow \infty }q_{k}=q.
\end{equation}

{\rm (a)} Then uniformly in compact subsets of $\{z:|z|<\frac{1}{%
2+|1+q|}\}${\rm ,}
\begin{equation}
\lim_{k\rightarrow \infty }H_{q_{k}}(z)=H_{q}(z).
\end{equation}

{\rm (b)} Let $\ell \geq 1$  and let $q$  be a primitive 
 $\ell^{\rm th}$ root of unity{\rm ,} and   
\begin{equation}
\rho (q_{k})>2^{-2/\ell },\qquad k\geq 1\mathit{.}
\end{equation}
  Let  $\Omega _{1}$  and $\Omega _{2}$  be open
connected sets with $\Omega _{1}\subseteq \Omega _{2}$  and $%
\Omega _{1}$ containing a branchpoint of  $H_{q}${\rm ,} that is{\rm ,}
containing one of the $\ell $ values of  $(-\frac{1}{4})^{1/\ell }$. Assume moreover that 
\begin{equation}
z\in \Omega _{1}\Rightarrow zq^{\pm 1}\in \Omega _{2}\mathit{.}
\end{equation}
 Then for large enough  $k${\rm ,} $H_{q_{k}}$  has a
pole in $\Omega _{2}$. If  $\ell $ is odd{\rm ,} and  $%
1>r>2^{-2/\ell }$  and  $\delta >0${\rm ,} then for large $k${\rm ,} $H_{q_{k}}$  has a pole in  $\{z:r<\left| z\right|
<r+\delta \}$.
\endproclaim

Thus (b) shows that every branchpoint of $H_{q}$ attracts a growing number
of poles of $H_{q_{k}}$ as $k\rightarrow \infty $. Next, we turn to
convergence of the c.f. Let us recall that we denoted the $n^{\rm th}$ convergent
for (1.6) by 
$$
\frac{\mu _{n}(z)}{\nu _{n}(z)}=1+\frac{qz|  }{|  1}+\frac{q^{2}z|  }{%
|  1}+\cdots +\frac{q^{n}z|  }{|  1}.
$$
It is known \cite{Hirschhorn1972} that 
\begin{equation}
\mu _{n}(z)=\sum_{k=0}^{[\frac{n+1}{2}]}z^{k}q^{k^{2}}\left[ 
\begin{array}{c}
n+1-k \\ 
k
\end{array}
\right]
\end{equation}
and 
\begin{equation}
\nu _{n}(z)=\mu _{n-1}(qz)=\sum_{k=0}^{[\frac{n}{2}]}z^{k}q^{k(k+1)}\left[ 
\begin{array}{c}
n-k \\ 
k
\end{array}
\right]
\end{equation}
where $[x]$ is the greatest integer $\leq x$, and 
$$
\left[ 
\begin{array}{c}
\alpha \\ 
l
\end{array}
\right] =\frac{(1-q^{\alpha })(1-q^{\alpha -1})\cdots (1-q^{\alpha -l+1})}{%
(1-q)(1-q^{2})\cdots (1-q^{l})},\qquad l\geq 0,\alpha \in {\Bbb C},
$$
is the Gaussian binomial coefficient. We shall reproduce the elegant proof
due to Adiga et al.\ \cite{Adigaetal1985} in Section 2.

When $G_{q}$ has positive radius of convergence, subsequences of the
numerators $\left\{ \mu _{n}\right\} _{n=1}^{\infty }$ and denominators $%
\left\{ \nu _{n}\right\} _{n=1}^{\infty }$ of the continued fraction
converge separately, depending on the behaviour of $q^{n}$. Of course, if $q$
is not a root of unity, then $\left\{ q^{n}\right\} _{n=1}^{\infty }$ is
dense on the unit circle, and one may extract a subsequence converging to an
arbitrary $\beta $ on the unit circle.

\nonumproclaim{Theorem 2.4} 
 Let  $q=e^{2\pi i\tau }${\rm ,} $\tau $ irrational. Let 
$|\beta |=1$ and  ${\cal S}$ be any infinite
sequence of positive integers with 
\begin{equation}
\lim_{n\rightarrow \infty ,n\in {\cal S}}q^{n}=\beta .
\end{equation}
 Then uniformly in compact subsets of $\{z:|z|<R(q)\}${\rm ,}
\begin{equation}
\lim_{n\rightarrow \infty ,n\in {\cal S}}\mu _{n}(z)=\overline{G_{q}(%
\overline{\beta qz})}G_{q}(z);
\end{equation}
\begin{equation}
\lim_{n\rightarrow \infty ,n\in {\cal S}}\nu _{n}(z)=\overline{G_{q}(%
\overline{\beta qz})}G_{q}(qz).
\end{equation}
 Moreover{\rm ,} uniformly in compact subsets of  $\{z:|z|<R(q)\}$ 
  omitting zeros of $G_{q}(\overline{\beta qz})$  and  $%
G_{q}(qz),$%
\begin{equation}
\lim_{n\rightarrow \infty ,n\in {\cal S}}\frac{H_{q}(z)-\frac{\mu _{n}(z)%
}{\nu _{n}(z)}}{(-1)^{n}z^{n+1}q^{(n+1)(n+2)/2}}=\frac{G_{q}(\beta q^{2}z)}{%
G_{q}(qz)^{2}\overline{G_{q}(\overline{\beta qz})}}
\end{equation}
 and so in such sets omitting these zeros{\rm ,} 
\begin{equation}
\lim_{n\rightarrow \infty ,n\in {\cal S}}\frac{\mu _{n}(z)}{\nu _{n}(z)}%
=H_{q}(z).
\end{equation}
\endproclaim

The crucial point in the last line is that the convergence takes place away
from the zeros of both $G_{q}(z)$ and $G_{q}(\overline{\beta qz})$. The
zeros of $G_{q}(\overline{\beta qz})$ need not be poles of $H_{q}$, and yet
(2.16) shows that they attract poles of the convergents. Moreover as the
zeros of $H_q$ are simply rotations and reflections of the zeros of $G_{q}(z)$ it
follows that if $H_{q}$ has poles of total multiplicity $\ell $ on a given
circle, then for all large enough $n$, $\frac{\mu _{n}}{\nu _{n}}$ has $%
2\ell $ poles close to this circle, that is, twice as many poles as the
approximated function! We formalize this as: 

\nonumproclaim{{C}orollary 2.5} Let $q=e^{2\pi i\tau },\tau $ be irrational. Assume that  $%
r<R(q)$ and  $H_{q}$ has poles of total multiplicity  $%
\ell $  on $\{z:|z|=r\}$. Let  $U$  be an
open set containing this circle. Then there exists $n_{0}$  such
that for  $n\geq n_{0}${\rm ,} $\mu _{n}/\nu _{n}$ has poles of
total multiplicity $\geq 2\ell $ in $U.$\endproclaim

This is the first such example in the literature, in which \textit{all}
approximants of large order have more poles than the approximated function
in a region of meromorphy. If we could show that there does not exist $\beta 
$ for which the zero sets of $G_{q}(qz)$ and $G_{q}(\overline{\beta qz})$
are the same, this would establish a counterexample to the Baker-Gammel-Wills
conjecture. For then, given any subsequence of the convergents, we can
extract a further subsequence for which (2.14) holds for some $\beta $; that
subsequence cannot converge uniformly in a compact set containing zeros of $%
G_{q}\left( \overline{\beta qz}\right) $ that are not zeros of $G_{q}\left(
z\right) $. For special $q$, we shall do this in Section 7.

Another feature of Theorem 2.4 is that it describes what happens only in $%
|z|<R(q)$, yet the function $H_{q}$ may have meaning in a much
larger circle. If for example $R(q)<\frac{1}{4}$, then $G_{q}$ is not
defined in $R(q)<|z|<\frac{1}{4}$, but by Worpitzky, $H_{q}$ is
analytic in $|z|<\frac{1}{4}$. One might hope for an alternative
formulation of (2.15) and (2.16) in this case. However this is not possible.
Those separate limits guarantee normality and uniform boundedness of $\{\mu
_{n}\}$ and $\{\nu _{n}\}$ in $|z|\leq r<R(q)$, but the following
result shows that $\{\mu _{n}\}$ and $\{\nu _{n}\}$ cannot be uniformly
bounded in $|z|\leq r$ for any $r>R(q)$. 

\nonumproclaim{Theorem 2.6} Let  $q=e^{2\pi i\tau }${\rm ,} $\tau $  irrational. Then
for $0<r<R(q),$%
\begin{equation}
\sup_{n\geq 1}\Vert \mu _{n}\Vert _{L_{\infty }(|z|\leq
r)}\ <\infty ;\ \sup_{n\geq 1}\Vert \nu _{n}\Vert _{L_{\infty }(|z|\leq r)}\ <\infty
\end{equation}
 and for $r>R(q),$%
\begin{equation}
\sup_{n\geq 1}\Vert \mu _{n}\Vert _{L_{\infty }(|z|\leq
r)}\ =\infty ;\sup_{n\geq 1}\Vert \nu _{n}\Vert _{L_{\infty }(|z|\leq r)}\ =\infty .
\end{equation}
\endproclaim

Thus in the case $R(q)<\frac{1}{2+|1+q|}$, the numerators $\left\{
\mu _{n}\right\} _{n=1}^{\infty }$ and denominators $\left\{ \nu
_{n}\right\} _{n=1}^{\infty }$ are normal in $\left\{ z:|z|<R(q)\right\} $, while in $\left\{ z:R(q)<|z|<\frac{1}{2+|1+q|}\right\} $, the numerators and denominators do not converge
separately, nor can they be normal, yet their ratio converges to $H_{q}$.

\section{Preliminaries}

In this section, we gather some elementary identities from the theory of
continued fractions, and also prove (2.12), (2.13) and some functional
relations for $G_{q}$. For the reader's convenience, we include many of the
proofs. Recall the notation $(a;q)_{0}:=1$ and 
\begin{equation}
(a;q)_{l}:=\prod_{j=1}^{\ell }(1-aq^{j-1}),\qquad \ell \geq 1.
\end{equation}

\nonumproclaim{Lemma 3.1} Let $\mu _{n}$  and $\nu _{n}$  be given by
{\rm (2.12), (2.13).} Then 
$$
\frac{\mu _{n}(z)}{\nu _{n}(z)}=1+\frac{qz|  }{|  1}+\frac{q^{2}z|  }{%
|  1}+\frac{q^{3}z|  }{|  1}+\cdots \frac{q^{n}z|  }{|  1}.
$$
\endproclaim

\demo{Proof} Fix $n\geq 1$. Following \cite[p.\ 26]{Adigaetal1985}, we set
for $r\geq 0,$%
$$
F_{r}:=\sum_{k=0}^{\infty }\frac{z^{k}q^{k(r+k)}(q;q)_{n-r-k+1}}{%
(q;q)_{k}(q;q)_{n-r-2k+1}}=\sum_{k=0}^{\infty }z^{k}q^{k(r+k)}\left[ 
\begin{array}{c}
n-r-k+1 \\ 
k
\end{array}
\right] .
$$
Then 
\begin{eqnarray*}
F_{r}-F_{r+1}&=&\sum_{k=0}^{\infty }\frac{z^{k}q^{k(r+k)}(q;q)_{n-r-k}}{%
(q;q)_{k}(q;q)_{n-r-2k}}\left[\frac{1-q^{n-r-k+1}}{1-q^{n-r-2k+1}}-q^{k}\right]
\\[4pt]
&=&\sum_{k=1}^{\infty }\frac{z^{k}q^{k(r+k)}(q;q)_{n-r-k}}{%
(q;q)_{k}(q;q)_{n-r-2k+1}}(1-q^{k})\\[4pt]
&=&\sum_{j=0}^{\infty }\frac{%
z^{j+1}q^{(j+1)(r+j+1)}(q;q)_{n-r-j-1}}{(q;q)_{j}(q;q)_{n-r-2j-1}}
 =zq^{r+1}F_{r+2},
\end{eqnarray*}
and hence 
$$
F_{r}/F_{r+1}=1+\frac{q^{r+1}z}{F_{r+1}/F_{r+2}}.
$$
Moreover,  we see easily  using $\left[ 
\begin{array}{c}
m \\ 
l
\end{array}
\right] =0,m>l$, that $F_{0}=\mu _{n};F_{1}=\nu _{n};F_{n-1}=1+zq^{n};F_{n}=1$. \pagebreak So 
\begin{eqnarray*}
\frac{\mu _{n}}{\nu _{n}}&=&F_{0}/F_{1}=1+\frac{qz}{F_{1}/F_{2}}
\\
&=&1+\frac{qz|  }{|  1}+\frac{q^{2}z|  }{|  F_{2}/F_{3}}=\cdots
\\
&=&1+\frac{qz|  }{|  1}+\frac{q^{2}z|  }{|  1}+\frac{q^{3}z|  }{| 
1}+\cdots \frac{q^{n}z|  }{|  1}.\\
\noalign{\vskip-36pt}
\end{eqnarray*}
\enddemo
 
\phantom{aman}

Next, we record the standard recurrence relations for the continued fraction
numerators and denominators:\newline

\nonumproclaim{Lemma 3.2} 
\begin{eqnarray}
\mu _{n}(z)&=&\mu _{n-1}(z)+q^{n}z\mu _{n-2}(z);
\\
\nu _{n}(z)&=&\nu _{n-1}(z)+q^{n}z\nu _{n-2}(z);
\\
(\mu _{n}\nu _{n-1}-\mu _{n-1}\nu _{n})(z)&=&(-1)^{n-1}z^{n}q^{\frac{n(n+1)}{2}%
}.
\end{eqnarray}
\endproclaim

\demo{Proof}  The first two are the standard recurrence relations for the
numerator and denominator of a continued fraction \cite[p.\ 20]{JonesThron1980}%
, \cite[pp.\ 8--9]{LorentzenWaadeland1992} though they may also be easily
proved from (2.12), (2.13) and the identity 
$$
\left[ 
\begin{array}{c}
m \\ 
l
\end{array}
\right] =\left[ 
\begin{array}{c}
m-1 \\ 
l
\end{array}
\right] +q^{m-l}\left[ 
\begin{array}{c}
m-1 \\ 
l-1
\end{array}
\right] .
$$
The third is also a standard relation, and is an easy consequence of (3.2),
(3.3).\enddemo

Next, we record an error formula for the difference between $H_{q}$ and the
convergents to its c.f., making use of the functional equation in the
process:

\nonumproclaim{Lemma 3.3} 
\begin{eqnarray}
\left(H_{q}-\frac{\mu _{n}}{\nu _{n}}\right)(z)&=&\frac{(-1)^{n}z^{n+1}q^{(n+1)(n+2)/2}}{%
\nu _{n}(z)[\nu _{n+1}(z)+\nu _{n}(z)[H_{q}(q^{n+1}z)-1]]}
\\
&=&\frac{(-1)^{n}z^{n+1}q^{(n+1)(n+2)/2}}{\nu _{n}(z)[\nu
_{n}(z)H_{q}(q^{n+1}z)+q^{n+1}z\nu _{n-1}(z)]}.
\end{eqnarray}
\endproclaim

\demo{Proof} We use the following elementary result in the theory of
continued fractions: let $\left\{ a_{j}\right\} ,\left\{ b_{j}\right\} $ be
complex numbers and 
$$
\frac{A_{k}}{B_{k}}=b_{0}+\frac{a_{1}|  }{|  b_{1}}+\cdots \frac{a_{k}|  }{%
|  b_{k}},\qquad k\geq 0.
$$
Then for $u\in {\Bbb C}$, 
\begin{equation}
 b_{0}+\frac{a_{1}|  }{|  b_{1}}+\cdots \frac{a_{k}|  }{|  b_{k}+u} 
 = b_{0}+\frac{a_{1}|  }{|  b_{1}}+\cdots \frac{a_{k}|  }{|  b_{k}}+%
\frac{u|  }{|  1}=\frac{A_{k}+A_{k-1}u}{B_{k}+B_{k-1}u}.\quad
\end{equation}
This follows immediately from the recurrence relations for the numerators and
denominators of the continued fraction. See for example \cite[p.\ 20]
{JonesThron1980}, \cite[p.\ 8]{LorentzenWaadeland1992}. Now in our situation,
our iterated functional relation (1.5) for $H_{q}$ gives 
$$
H_{q}\left( z\right) =1+\frac{qz|  }{|  1}+\frac{q^{2}z|  }{|  1}+\cdots %
\frac{q^{n}z|  }{|  1+u},
$$
where $u:=H_{q}(q^{n}z)-1$. By (3.7),
\begin{equation}
H_{q}(z)=\frac{\mu _{n}(z)+\mu _{n-1}(z)[H_{q}(q^{n}z)-1]}{\nu _{n}(z)+\nu
_{n-1}(z)[H_{q}(q^{n}z)-1]}.
\end{equation}
Then 
\begin{eqnarray*}
\left(H_{q}-\frac{\mu _{n-1}}{\nu _{n-1}}\right)(z)&=&\frac{(\mu _{n}\nu _{n-1}-\mu
_{n-1}\nu _{n})(z)}{\nu _{n-1}(z)[\nu _{n}(z)+\nu _{n-1}(z)[H_{q}(q^{n}z)-1]]%
}
\\
&=&\frac{(-1)^{n-1}z^{n}q^{n(n+1)/2}}{\nu _{n-1}(z)[\nu _{n}(z)+\nu
_{n-1}(z)[H_{q}(q^{n}z)-1]]}.
\end{eqnarray*}
Replacing $n-1$ by $n$ gives the first identity (3.5) and then our
recurrence relation (3.3) gives the second relation (3.6).\enddemo

Next, we establish some functional relations for $G_{q}$:
\nonumproclaim{Lemma 3.4} Let  $q=e^{2\pi i\tau }${\rm ,} with $\tau $  irrational. Then 
\begin{equation}
G_{q}(z)=G_{q}(qz)+qzG_{q}(q^{2}z).
\end{equation}
 Moreover if  $\ell \geq 1${\rm ,} 
\begin{equation}
G_{q}\left(\frac{z}{q^{\ell }}\right)=G_{q}(qz)\mu _{\ell }\left(\frac{z}{q^{\ell }}%
\right)+qzG_{q}(q^{2}z)\mu _{\ell -1}\left(\frac{z}{q^{\ell }}\right).
\end{equation}
\endproclaim

\demo{Proof} Firstly, (3.9) follows easily from the series definition of $%
G_{q}$. We prove (3.10) by induction on $\ell $. If  we define $\mu _{-1}:=1$, then (3.10) follows from (3.9) for $\ell =0$.
Assume now as an induction hypothesis that (3.10) is true for $\ell $. Then using our recurrence
relation for $\mu _{\ell +1}$, 
\begin{eqnarray*}
&&G_{q}(qz)\mu _{\ell +1}\left(\frac{z}{q^{\ell +1}}\right)+qzG_{q}(q^{2}z)\mu _{\ell }(%
\frac{z}{q^{\ell +1}})
\\
&&\qquad =G_{q}(qz)\left[\mu _{\ell }\left(\frac{z}{q^{\ell +1}}\right)+z\mu _{\ell -1}\left(\frac{z}{%
q^{\ell +1}}\right)\right]+qzG_{q}(q^{2}z)\mu _{\ell }\left(\frac{z}{q^{\ell +1}}\right)
\\
&&\qquad =G_{q}(z)\mu _{\ell }\left(\frac{z}{q^{\ell +1}}\right)+zG_{q}(qz)\mu _{\ell -1}\left(\frac{z%
}{q^{\ell +1}})=G_{q}(\frac{z}{q^{\ell +1}}\right)
\end{eqnarray*}
by first (3.9) and then our induction hypothesis that (3.10) is true for $%
\ell $. So we have the result for $\ell +1$.\enddemo

Note that if we define $H_{q}$ by (1.3), then the functional relation (1.4)
follows immediately from (3.9).

\section{Roots of unity}

In this section we prove Theorem 2.1 using the following result 
\cite[Satz 2.38, p.~86]{Perron1957}:

\nonumproclaim{Lemma 4.1} Consider the {\rm c.f.}\
\begin{equation}
b_{0}+\frac{a_{1}|  }{|  b_{1}}+\cdots \frac{a_{l-1}|  }{|  b_{l-1}}+%
\frac{a_{l}|  }{|  b_{0}}+\frac{a_{1}|  }{|  b_{1}}+\cdots \frac{%
a_{l-1}|  }{|  b_{l-1}}+\frac{a_{l}|  }{|  b_{0}}+\frac{a_{1}|  }{%
|  b_{1}}+\cdots ,
\end{equation}
 periodic of period  $\ell $\textit{. Let }$A_{k}/B_{k}$ 
denote the $k^{\rm th}$ convergent{\rm,} $k\geq 0$.  Let  $B_{\ell
-1}\neq 0${\rm ,} and  $x_{1},x_{2}$ denote the roots of the
quadratic 
\begin{equation}
B_{\ell -1}x^{2}+(a_{\ell }B_{\ell -2}-A_{\ell -1})x-a_{\ell }A_{\ell -2}=0.
\end{equation}
 Assume either that {\rm (a)} $x_{1}=x_{2}$  or {\rm (b)} $x_{1}\neq
x_{2}$  and both the following hold\/{\rm :} 
\begin{equation}
|B_{\ell -1}x_{1}+a_{\ell }B_{\ell -2}| >  |B_{\ell
-1}x_{2}+a_{\ell }B_{\ell -2}|  ;
\end{equation}
\begin{equation}
A_{k}-x_{2}B_{k}\neq 0,\qquad k=0,1,2,\ldots \ell -2.
\end{equation}
 Then 
$$
\lim_{k\rightarrow \infty }A_{k}/B_{k}=x_{1}.
$$
\endproclaim

Of course in our case $a_{j}=q^{j}z$; $b_{j}=1$; $A_{j}=\mu _{j}$; $%
B_{j}=\nu _{j}$. Before applying the above result, we need
\nonumproclaim{Lemma 4.2} Assume that $q$ is a primitive $\ell^{\rm th}$ root of
unity. For our choice of  $\{a_{j}\},\{b_{j}\}${\rm ,} 
\begin{equation}
a_{\ell }B_{\ell -2}+A_{\ell -1}=z\nu _{\ell -2}(z)+\mu _{\ell -1}(z)=1.
\end{equation}
\endproclaim

\demo{Proof} We shall use the explicit forms (2.12), (2.13) for $\mu
_{n},\nu _{n}$. We have 
\begin{eqnarray*}
z\nu _{\ell -2}(z)+\mu _{\ell -1}(z)&=&\sum_{k=0}^{\left[ \frac{\ell -2}{2}%
\right] }z^{k+1}q^{k(k+1)}\left[ 
\begin{array}{c}
\ell -2-k \\ 
k
\end{array}
\right] +\sum_{k=0}^{\left[ \frac{\ell }{2}\right] }z^{k}q^{k^{2}}\left[ 
\begin{array}{c}
\ell -k \\ 
k
\end{array}
\right]
\\
&=&1+\sum_{j=1}^{\left[ \frac{\ell }{2}\right] }z^{j}q^{j(j-1)}\left( \left[ 
\begin{array}{c}
\ell -1-j \\ 
j-1
\end{array}
\right] +q^{j}\left[ 
\begin{array}{c}
\ell -j \\ 
j
\end{array}
\right] \right) .
\end{eqnarray*}
Now we see that 
\begin{eqnarray*}
\left[ 
\begin{array}{c}
\ell -1-j \\ 
j-1
\end{array}
\right] +q^{j}\left[ 
\begin{array}{c}
\ell -j \\ 
j
\end{array}
\right] &=&\left[ 
\begin{array}{c}
\ell -1-j \\ 
j-1
\end{array}
\right] (1+\frac{q^{j}(1-q^{\ell -j})}{1-q^{j}})
\\
&=&\left[ 
\begin{array}{c}
\ell -1-j \\ 
j-1
\end{array}
\right] \frac{1-q^{\ell }}{1-q^{j}}=0
\end{eqnarray*}
as $q$ is a primitive $\ell ^{\rm th}$ root of unity.\enddemo

We turn to the

\demo{Proof of Theorem {\rm 2.1}} The quadratic (4.2) becomes (recall $q^{\ell
}=1$) 
\begin{equation}
\nu _{\ell -1}x^{2}+(z\nu _{\ell -2}(z)-\mu _{\ell -1}(z))x-z\mu _{\ell
-2}(z)=0.
\end{equation}
The discriminant of this quadratic is 
\begin{eqnarray*}
D&:=&(z\nu _{\ell -2}(z)-\mu _{\ell -1}(z))^{2}+4z\mu _{\ell -2}(z)\nu _{\ell
-1}(z)
\\
&=&([z\nu _{\ell -2}(z)+\mu _{\ell -1}(z)]-2\mu _{\ell -1}(z))^{2}+4z\mu
_{\ell -2}(z)\nu _{\ell -1}(z)
\\
&=&[z\nu _{\ell -2}(z)+\mu _{\ell -1}(z)]^{2}-4[z\nu _{\ell -2}(z)+\mu _{\ell
-1}(z)]\mu _{\ell -1}(z)\\
&&+\ 4\mu _{\ell -1}(z)^{2}+4z\mu _{\ell -2}(z)\nu
_{\ell -1}(z)
\\
&=&1-4z(\mu _{\ell -1}\nu _{\ell -2}-\mu _{\ell -2}\nu _{\ell
-1})(z)=1-4z^{\ell }(-1)^{\ell -2}q^{\ell (\ell -1)/2},
\end{eqnarray*}
by first Lemma 4.2 and then (3.4). Next, we note that 
\begin{equation}
(-1)^{\ell -1}q^{\ell (\ell -1)/2}=1.
\end{equation}
Indeed for $\ell $ even, $q^{\ell /2}=-1$ (as $q$ is a primitive $\ell ^{\rm th}$
root of unity), and for $\ell $ odd, $(\ell -1)/2$ is an integer. Thus 
$$
D=1+4z^{\ell }
$$
and the roots of the quadratic (4.6) are, by Lemma 4.2, 
$$
x_{1}=\frac{-(z\nu _{\ell -2}-\mu _{\ell -1})(z)+\sqrt{D}}{2\nu _{\ell -1}(z)%
}=\frac{\mu _{\ell -1}-\frac{1}{2}+\sqrt{D/4}}{\nu _{\ell -1}};
$$
$$
x_{2}=\frac{-(z\nu _{\ell -2}-\mu _{\ell -1})(z)-\sqrt{D}}{2\nu _{\ell -1}(z)%
}=\frac{\mu _{\ell -1}-\frac{1}{2}-\sqrt{D/4}}{\nu _{\ell -1}}.
$$
(The branch of the $\sqrt{}$ is the principal one.) Now we examine (4.3) and
(4.4). Firstly (4.3) becomes 
$$
|\nu _{\ell -1}x_{1}+z\nu _{l-2}| >| \nu _{\ell -1}x_{2}+z\nu
_{l-2}|,
$$
that is, in view of Lemma 4.2, 
$$
\left| \frac{1}{2}+\sqrt{D/4}\right| >\left| \frac{1}{2}-\sqrt{D/4}\right| .
$$
Now by our choice of the principal branch of $\sqrt{}$, $\sqrt{D/4}=\alpha
+i\beta $, where $\alpha >0$, provided $D\notin (-\infty ,0]$. Then 
$$
\left|  \frac{1}{2}+\sqrt{D/4}\right| ^{2}=\left(\frac{1}{2}+\alpha \right)^{2}+\beta ^{2}>\left( 
\frac{1}{2}-\alpha \right)^{2}+\beta ^{2}=\left| \frac{1}{2}-\sqrt{D/4}\right| ^{2}.
$$
So we have (4.3) provided 
$$
D=1+4z^{\ell }\notin (-\infty ,0]\Leftrightarrow z^{\ell }\notin \left(-\infty ,-%
\frac{1}{4}\right]\Leftrightarrow z\notin {\cal L}.
$$
Next, we examine (4.4). We see that for $\nu _{\ell -1}(z)\neq 0,$ 
\begin{eqnarray}
A_{k}-x_{2}B_{k}\neq 0&\Leftrightarrow& \mu _{k}-\frac{\nu _{k}}{\nu _{\ell -1}%
}\left[\mu _{\ell -1}-\frac{1}{2}-\sqrt{D/4}\right]\neq 0
\\
&\Leftarrow& \left(\mu _{k}\nu _{\ell -1}-\left[\mu _{\ell -1}-\frac{1}{2}\right]\nu _{k}\right)^{2}-%
\frac{D}{4}\nu _{k}^{2}\neq 0
\nonumber\\
&\Leftrightarrow &J_{k}^{2}+J_{k}\nu _{k}-z^{\ell }\nu _{k}^{2}\neq 0,\nonumber
\end{eqnarray}
where 
$$
J_{k}:=\mu _{k}\nu _{\ell -1}-\mu _{\ell -1}\nu _{k}.
$$
Now (3.4) shows that for each $n$, $\frac{\mu _{n}}{\nu _{n}}-\frac{\mu
_{n-1}}{\nu _{n-1}}$ has a zero of order $n$ at $0$. Adding this for $%
n=k+1,k+2,\ldots,\ell -1$, shows that 
$$
\frac{J_{k}}{\nu _{k}\nu _{\ell -1}}=\frac{\mu _{k}}{\nu _{k}}-\frac{\mu
_{\ell -1}}{\nu _{\ell -1}}
$$
has a zero of order at least $k+1$ at $0$. Thus 
$$
J_{k}=J_{k}\left( z\right) =z^{k+1}\pi _{k}\left( z\right) ,
$$
where, as $\deg (\mu _{k})\leq \frac{k+1}{2};\deg (\nu _{k})\leq \frac{k}{2}$%
, we see that $\pi _{k}$ is a polynomial of degree at most $\frac{k+\ell }{2}%
-\left( k+1\right) =\frac{\ell -2-k}{2}$. So (4.8)  becomes, after division
by $z^{k+1}$, $\phantom{\sum^\int}$
\begin{equation}
z^{k+1}\pi _{k}^{2}+\pi _{k}\nu _{k}-z^{\ell -k-1}\nu _{k}^{2}\neq 0.
\end{equation}
(Recall that the c.f.\ converges at $z=0$, so that point can be omitted.) The
left-hand side of (4.9) is a polynomial of degree $\leq \ell -1$, so has at
most $\ell -1$ zeros. Considering $k=0,1,2,\ldots,\ell -2$, we obtain a set $%
{\cal P}$ of at most $(\ell -1)(\ell -1)$ exceptional points. Adding the
at most $(\ell -1)/2$ zeros of $\nu _{\ell -1}$  gives a set ${\cal P}$
of at most $\frac{(\ell -1)(2\ell -1)}{2}$ points.\enddemo

\section{Proof of Theorems 2.2, 2.3}

The proof of Theorem 2.2 requires three lemmas. The first is a special case
of a theorem of P\'olya: 

\nonumproclaim{Lemma 5.1} 
 Let  $g$  be a function meromorphic in $| z| <\sigma ${\rm ,} and let 
$g$ be analytic at $0${\rm ,} with
Maclaurin series $\sum_{j=0}^{\infty }g_{j}z^{j}$. Let 
\begin{equation}
D_{n}(g):=\det (g_{1+i+j})_{i,j=0}^{n-1}.
\end{equation}
 Then  
\begin{equation}
\limsup_{n\rightarrow \infty }| D_{n}(g)| ^{1/n^{2}}\leq \sigma ^{-1}.
\end{equation}
\endproclaim

\demo{Proof} This first appeared in \cite{Polya1928}. A more accessible
reference is \cite[p.~305, Thm.~3]{Goluzin1969}, though the proof there is
for analytic $f$. \enddemo

Next, we record the well-known relation between the Hankel determinants $%
D_{n}(g)$ and the continued fraction coefficients of $g$. It was used for
example in~\cite{ArmsEdrei1970}:

\nonumproclaim{Lemma 5.2} Assume that $g$ is analytic near  $0${\rm ,} and has
{\rm (}\/formal\/{\rm )} continued fraction expansion 
$$
c_{0}+\frac{c_{1}z|  }{|  1}+\frac{c_{2}z|  }{|  1}+\frac{c_{3}z|  
}{|  1}+\cdots 
$$
 with all $c_{j}\neq 0$. Then 
\begin{equation}
D_{n}(g)=c_{1}^{n}\prod_{j=1}^{n-1}(c_{2j}c_{2j+1})^{n-j}.
\end{equation}
\endproclaim

\demo{Proof} If $g$ has Maclaurin series coefficients $\{g_{j}\}$ and we
define 
$$
H_{k}^{(\ell )}:=\det (g_{\ell +i+j})_{i,j=0}^{k-1}
$$
then is it known \cite[p.\ 257]{LorentzenWaadeland1992} that 
$$
c_{2k}=-\frac{H_{k-1}^{(1)}H_{k}^{(2)}}{H_{k}^{(1)}H_{k-1}^{(2)}};c_{2k+1}=-%
\frac{H_{k+1}^{(1)}H_{k-1}^{(2)}}{H_{k}^{(1)}H_{k}^{(2)}}.
$$
We deduce that 
$$
\prod_{k=1}^{\ell }(c_{2k}c_{2k+1})=\prod_{k=1}^{\ell }\frac{%
H_{k-1}^{(1)}H_{k+1}^{(1)}}{(H_{k}^{(1)})^{2}}=\frac{H_{\ell +1}^{(1)}}{%
H_{\ell }^{(1)}}/\frac{H_{1}^{(1)}}{H_{0}^{(1)}}=\frac{H_{\ell +1}^{(1)}}{%
H_{l}^{(1)}}/c_{1}.
$$
Multiplying this for $\ell =1,2,\ldots,n-1$ and noting that $%
D_{n}(g)=H_{n}^{(1)} $ and $H_{1}^{(1)}=c_{1}$ gives the result.\enddemo

Next, we record one form of the \textit{fundamental inequalities}, as a
criterion for convergence of continued fractions:
\nonumproclaim{Lemma 5.3} Assume that the {\rm c.f. } 
$$
\frac{1|}{|  1}+\frac{a_{2}|}{|  1}+\frac{a_{3}|}{|  1}+\frac{a_{4}|}{%
|  1}+\cdots 
$$
 satisfies for some sequence  $\left\{ r_{j}\right\} _{j=1}^{\infty
}\subset \left( 0,\infty \right) $ the fundamental inequalities  
\begin{equation}
r_{j}| 1+a_{j}+a_{j+1}| \geq r_{j}r_{j-2}| a_{j}| +|
a_{j+1}| ,\qquad j\geq 1
\end{equation}
 where  
$$
a_{1}:=0;\quad r_{0}:=0;\quad r_{-1}:=0.
$$
 Let  $A_{j}/B_{j}$\textit{\ denote the }$j^{\rm th}$ convergent to
the {\rm c.f.}\ above. Then $B_{j}\neq 0,j\geq 1${\rm ,} and 
\begin{equation}
\left| \frac{A_{j+1}}{B_{j+1}}-\frac{A_{j}}{B_{j}}\right| \leq
r_{1}r_{2}\cdots r_{j},\qquad j\geq 1.
\end{equation}
\endproclaim

\demo{Proof} This is Theorem 9.1 in \cite[p.\ 41]{Wall1973} and inequality
(9.4) in \cite[p.~42]{Wall1973}. \phantom{soonIwill}\enddemo

Now we turn to the

\demo{Proof of Theorem {\rm 2.2(a)}} We first show that the c.f.\ (1.6) converges
to a function $H^{*}(z)$ analytic in $\{z:| z| <\frac{1}{2+|
1+q| }\}$. This part works even if $q$ is a root of unity.  Let $K$ be a compact subset of this ball. Choose $\varepsilon
>0$ such that 
$$
\left| z\right| <\frac{1-\varepsilon }{2+\left| 1+q\right| },\qquad z\in K.
$$
We apply the fundamental inequalities (5.4) with 
\begin{eqnarray*}
a_{j} &:=&q^{j}z,\qquad j\geq 2; \\
r_{j} &:=&1-\varepsilon ,\qquad j\geq 1.
\end{eqnarray*}
For $j=1$, we see that 
\begin{eqnarray*}
r_{j}  \left|   1+a_{j}+a_{j+1}\right| &=&\left( 1-\varepsilon \right) \left|
1+q^{2}z\right| \\
&\geq &\frac{1-\varepsilon }{2}>\left| z\right| =r_{1}r_{-1}\left|
a_{1}\right| +\left| a_{2}\right| .
\end{eqnarray*}
For $j\geq 2$, 
\begin{eqnarray*}
r_{j}  \left|   1+a_{j}+a_{j+1}\right|& =&\left( 1-\varepsilon \right) \left|
1+q^{j}z\left( 1+q\right) \right| \\
&\geq &\left( 1-\varepsilon \right) \left( 1-\frac{\left| 1+q\right| }{%
2+\left| 1+q\right| }\right) =\frac{2\left( 1-\varepsilon \right) }{2+\left|
1+q\right| } \\
&>&2\left| z\right| \geq r_{j}r_{j-2}| a_{j}| +| a_{j+1}| .
\end{eqnarray*}
Thus the fundamental inequalities are satisfied. If $A_{j}\left( z\right)
/B_{j}\left( z\right) $ denotes the $j^{\rm th}$ convergent to the c.f.\ 
$$
\frac{1|}{|  1}+\frac{q^{2}z|}{|  1}+\frac{q^{3}z|}{|  1}+\frac{q^{4}z|%
}{|  1}+\cdots 
$$
then Lemma 5.3 shows that $B_{j}\left( z\right) \neq 0$ for $j\geq 1$ and $%
z\in K$, and 
$$
\left| \frac{A_{j+1}\left( z\right) }{B_{j+1}\left( z\right) }-\frac{%
A_{j}\left( z\right) }{B_{j}\left( z\right) }\right| \leq \left(
1-\varepsilon \right) ^{j},\qquad j\geq 1.
$$
Then $\left\{ A_{j}/B_{j}\right\} _{j=1}^{\infty }$ converges uniformly in $%
K $ and so the limit function is analytic in the interior of $K$. The same
is then true for the c.f.\ $H^{*}$ defined by (1.6). Thus $H^{*}$ is analytic
in $\{z:| z| <\frac{1}{2+| 1+q| }\}$, so 
$$
\rho (q)\geq \frac{1}{2+| 1+q| }.
$$
Next, we note that if $R(q)>0$, the function $H_{q}(z):=G_{q}(z)/G_{q}(qz)$
satisfies the same functional equation as does $H^{*}$, in view of the
functional equation~(3.9) for $G_{q}$. Moreover, $H^{*}(0)=H_{q}\left(
0\right) =1$. Then $H_{q}$ and $H^{*}$ have the same c.f.\ expansion and
hence the same Maclaurin series. This follows as the c.f.\ uniquely
determines the corresponding Maclaurin series. Then $H_{q}\left( z\right)
=G_{q}\left( z\right) /G_{q}\left( qz\right) $ provides a meromorphic
continuation of $H^{*}$ to $\left\{ z:\left| z\right| <R\left( q\right)
\right\} $. So we have the inequality 
$$
\rho (q)\geq \max \left\{\frac{1}{2+| 1+q| },R(q)\right\}.
$$
To show $\rho (q)\leq 1$, we note from Lemma 5.2 that 
$$
| D_{n}(H_{q})| =1,n\geq 1
$$
and from Lemma 5.1, 
$$
1=\limsup_{n\rightarrow \infty }| D_{n}(H_{q})| ^{1/n^{2}}\leq \frac{1%
}{\rho (q)}.
$$
Thus, we have $\rho (q)\leq 1$ and (2.7). To show that $H_{q}$ has a natural
boundary on $\{z:| z| =\rho (q)\}$, let us suppose that $z_{0}$ is a
point of analyticity of $H_{q}$ with $| z_{0}| =\rho (q)$. Then we can
find a ball $U$ centre $z_{0}$ in which $H_{q}$ is analytic and hence
meromorphic. The functional equation (1.4) in the form 
$$
H_{q}(qz)=\frac{qz}{H_{q}(z)-1}
$$
shows that $H_{q}(z)$ has a meromorphic continuation to the ball $%
qU=\break \{qz: q\in U\}$. Iteration of this argument shows that $H_{q}$ has a
meromorphic continuation to $q^{j}U=\{q^{j}z:z\in U\},j\geq 1$. As finitely
many such balls cover the circle $\left\{ z:| z| =\rho (q)\right\} $,
we obtain a meromorphic continuation of $H_{q}$ to $\{z:| z| <\rho
(q)+\varepsilon \}$, for some $\varepsilon >0$, contradicting the definition of $%
\rho (q)$. So $H_{q}$ must have a natural boundary on the circle $\{z:|
z| =\rho (q)\}$. \enddemo

In the proof of Theorem 2.2(b), we need part of a result of G.\ Petruska:%
\nonumproclaim{Lemma 5.4} Let  $c\in [0,\infty ]$. There exists an irrational number $%
\tau $  with continued fraction 
$$
\tau =\frac{1|  }{|  a_{1}}+\frac{1|  }{|  a_{2}}+\frac{1|  }{| 
a_{3}}+\cdots 
$$
 {\rm (}\/all $a_{i}$ positive integers\/{\rm )} such that if $\pi
_{n}/\chi _{n}$  denotes the  $n^{\rm th}$ convergent to the
{\rm c.f.}\
of  $\tau ${\rm ,} then 
\begin{equation}
\lim_{n\rightarrow \infty }\frac{\log \chi _{n+1}}{\chi _{n}}=c.
\end{equation}
\endproclaim

\demo{Proof} For the case $0<c<\infty $ this is part of Lemma 2 in 
\cite[p.~354]{Petruska1992} and for the case $c=\infty $, this was noted in 
\cite[p.~474, eqn. (1.17)]{Driveretal1991}. Almost every $\tau \in [0,1]$
satisfies (5.6) with $c=0$. Indeed this follows from Khinchin's theorem \cite
{Khinchin1936} that for a.e.\ $\tau $, 
\vglue12pt
\hfill ${\displaystyle
\lim_{n\rightarrow \infty }\frac{\log \chi _{n}}{n}=\frac{\pi ^{2}}{12\log 2}%
.}$
\enddemo
\vglue12pt

\demo{Proof of Theorem {\rm 2.2(b)}} Let us suppose that $G_{q}$ is analytic at
a point $z_{0}$ on the circle $|  z|  =R(q)$ and hence in some ball $U$
centre $z_{0}$. We shall use the functional relation (3.10) in the form 
$$
G_{q}\left(\frac{u}{q^{\ell +1}}\right)=G_{q}(u)\mu _{\ell }\left(\frac{u}{q^{\ell +1}}%
\right)+G_{q}(qu)u\mu _{\ell -1}\left(\frac{u}{q^{\ell +1}}\right).
$$
Let $1>\varepsilon >0$. Let us choose $\ell $ large with $\frac{z_{0}}{%
q^{\ell +1}}\in U$ and in fact such that the ball $U_{1}\;$centre $\frac{%
z_{0}}{q^{\ell +1}}$ and $1-\varepsilon $ times the radius of $U$ lies
inside $U$. Then the above identity shows that $G_{q}(qu)$ is meromorphic in 
$U_{1}$ and hence $G_{q}$ is meromorphic in $qU_{1}=\{qz:z\in U_{1}\}$.
Since $\varepsilon >0$ is arbitrary, we deduce that $G_{q}$ is meromorphic
in $qU$. By the same argument $G_{q}(u)$ is meromorphic in $q^{j}U,\ j\geq~1$. Finitely many such neighbourhoods cover
the circle $\left\{ z:| z| =R(q)\right\} $. Then $G_{q}$ is analytic on this circle, except possibly
for finitely many poles. Since there are only finitely many such poles, we
can choose $z_{0}$ with $| z_{0}| =R(q)$ such that both $z_{0}$ and $%
qz_{0}$ are points of analyticity. Thus there exists an open ball $B$ centre 
$z_{0}$, such that $G_{q}$ is analytic in both $B$ and $qB$. But the
functional relation (3.9) shows that $G_{q}$ is analytic in $q^{2}B$, and
hence also in $q^{j}B,j\geq 1.$ Hence $G_{q}$ is analytic on the whole
circle $\left\{ z:| z| =R(q)\right\} $, contradicting the definition
of $R(q)$. Thus, $G_{q}$ has a natural boundary on its circle of convergence.

Next we show that $R(q)$ may assume any value in $[0,1]$ by Petruska's
Lemma 5.4 with $q=e^{2\pi i\tau }$. We shall show (recall (2.5)) that 
$$
R(q)=\liminf_{n\rightarrow \infty }\Vert j\tau \Vert ^{1/j}=e^{-c}
$$
and since $e^{-c}$ may assume any value in $[0,1]$, the result follows. We
recall some elementary facts from the theory of continued fraction
expansions of real numbers: firstly if $\frac{j}{k}$ is not a convergent to
the c.f.\ of $\tau $, then \cite[p.\ 153, Thm.~184]{HardyWright1975} 
$$
\left| \tau -\frac{j}{k}\right| \geq \frac{1}{2k^{2}}
$$
and so if $k$ is not a denominator of some convergent, 
$$
\Vert k\tau \parallel\ \geq \frac{1}{2k}.
$$
Hence if we let ${\cal S}:=\{\chi _{1},\chi _{2},\chi _{3},\ldots\}$, then 
$$
\lim_{k\rightarrow \infty ,k\notin {\cal S}}\Vert k\tau \Vert ^{1/k}=1.
$$
Next for convergents $\frac{\pi _{n}}{\chi _{n}}$, we have the inequalities 
$$
\frac{1}{2\chi _{n+1}}\leq\ |  \chi _{n}\tau -\pi _{n}| <\frac{1}{\chi
_{n+1}}.
$$
See \cite[p.\ 140]{HardyWright1975} for the upper bound. The lower bound
follows from the error formula (see for example \cite[p.\ 354]{Petruska1992}) 
$$
\chi _{n}\tau -\pi _{n}=\frac{(-1)^{n}}{\chi _{n+1}+\alpha _{n+1}\chi _{n}},
$$
where $\alpha _{n+1}\in (0,1)$ and   $\chi _{n}$ increases with $n$. Then
we see that 
$$
\lim_{j\rightarrow \infty ,j\in {\cal S}}\Vert j\tau \Vert ^{1/j}\ =\lim_{n\rightarrow \infty }\Vert \chi _{n}\tau \Vert ^{1/\chi
_{n}}\ =\lim_{n\rightarrow \infty }\chi _{n+1}^{-1/\chi _{n}}=e^{-c},
$$
as desired.\enddemo

\demo{Proof of Theorem {\rm 2.2(c)}} 
When $R(q)=1$, of course $\rho (q)=1$ follows from (2.7), so that  Theorem 2.2(c)
follows immediately. Moreover, we noted in Section~2 that for a.e. $q$, $%
R\left( q\right) =1$.\enddemo

\demo{Proof of Theorem {\rm 2.3(a)}} We recall that we showed in the proof of
Theorem 2.2(a) that $H_{q}(z)\neq \infty ,|  z|  <\frac{1}{2+|  1+q|  }$ even if $q$ is a root of unity. Then our functional
equation (1.4), in the form 
\begin{equation}
\lbrack H_{q}(z)-1]H_{q}(qz)=qz
\end{equation}
shows that $H_{q}(qz)$ cannot have any zeros in that ball (recall that $%
H_{q}(0)\break =1$).  So $H_{q}$ does not assume the values $0,\infty $ there. By
the same token, the functional relation shows that $H_{q}$ cannot assume the
value $1$ in the punctured ball ${\cal B}:=\{z:0<|  z|  <\frac{1}{%
2+|  1+q|  }\}$ (for if $H_{q}(z)=1$, then $H_{q}(qz)=\infty $). Thus $%
H_{q}$ omits the values $0,1,\infty $ in that punctured ball. If $%
\{q_{k}\}_{k=1}^{\infty }$ satisfy (2.8), then in a given compact subset of $%
{\cal B}$, for large $k$, $H_{q_{k}}$ omits the values $0,1,\infty $ and
by Montel's theorem \cite[p.\ 54, p.\ 74]{Schiff1993}, $\{H_{q_{k}}\}_{k=1}^{%
\infty }$ is normal there. Let $H^{*}$ denote a limit function of some
subsequence, so that $H^{*}$ is either identically $\infty $ or is analytic
in ${\cal B}$. It follows easily from (5.7) that $H^{*}$ cannot be
identically $\infty $, so is analytic in ${\cal B}$. In view of (5.7) and
(2.8), we have the functional equation 
\begin{equation}
\lbrack H^{*}(z)-1]H^{*}(qz)=qz,z\in {\cal B}.
\end{equation}
Now for all $q_{k}$, the c.f.\ coefficients have absolute value $1$, so that by
Worpitzky's theorem \cite[p.\ 35]{LorentzenWaadeland1992} 
$$
|  H_{q_{k}}(z)-1| \leq \frac{1}{2},  |  z|   \leq \frac{1}{4}
$$
and so the same is true of $H^{*}$. It follows that $0$ is a removable
singularity of $H^{*}$ and defining $H^{*}(0)=1$, we obtain a function
analytic in $|  z|  <\frac{1}{2+|  1+q|  }$ satisfying the \pagebreak same
functional equation as $H_{q}$. Then both have the same c.f., both have the
same set of convergents, and hence the same Maclaurin series, so that  $%
H^{*}=H_{q} $. As $H^{*}$ was the limit of any subsequence, the full
sequence converges to $H_{q}$.\enddemo

\demo{Proof of Theorem {\rm 2.3(b)}} Let us assume that $\Omega _{2}$ does not
contain any pole of $H_{q_{k}}$ for infinitely many $k$. By passing to a
subsequence, we may assume that $H_{q_{k}}$ has no poles in $\Omega _{2}$
for all $k$. We may also assume that $z\in \Omega _{1}\Rightarrow
zq_{k}^{\pm 1}\in \Omega _{2}$ for all $k$ (if necessary make $\Omega _{1}$
a little smaller). Then our functional equation for $H_{q_{k}}$ in the form 
$$
\lbrack H_{q_{k}}(z)-1]H_{q_{k}}(q_{k}z)=q_{k}z
$$
shows that as $H_{q_{k}}(q_{k}z)\neq \infty ,z\in \Omega _{1}$, so that $%
H_{q_{k}}(z)\neq 1,z\in \Omega _{1}$. Similarly, 
$$
\lbrack H_{q_{k}}(z/q_{k})-1]H_{q_{k}}(z)=z
$$
and $H_{q_{k}}(z/q_{k})\neq \infty ,z\in \Omega _{1}$ implies $%
H_{q_{k}}(z)\neq 0,z\in \Omega _{1}$. Thus $H_{q_{k}}$ omit the values $%
0,1,\infty $ in $\Omega _{1}$ for all $k$, and so   $\left\{ H_{q_{k}}\right\}
_{k=1}^{\infty }$ is a normal family there. Let $H^{*}$ denote a
subsequential limit. As above it is not identically $\infty $, so is
analytic in $\Omega _{1}$. Then we have the functional equation (5.8) for $%
H^{*}$. Iterating the functional relation leads to the same continued
fraction as for $H_{q}$, periodic of period $\ell $. Moreover, the error
formula used in the proof of Lemma 3.3 gives (3.8) with $n=\ell $; that is, 
$$
H^{*}=\frac{\mu _{\ell }+\mu _{\ell -1}[H^{*}-1]}{\nu _{\ell }+\nu _{\ell
-1}[H^{*}-1]}.
$$
As before, $H^{*}(z)$ is one of the roots of the quadratic (4.6) (we also
use (3.2), (3.3)). Then as in the proof of Theorem 2.1, 
$$
H^{*}(z)=\frac{\mu _{\ell -1}(z)-\frac{1}{2}\pm \sqrt{\frac{1}{4}+z^{\ell }}%
}{\nu _{\ell -1}(z)}.
$$
It follows that we obtain a branch of $\sqrt{\frac{1}{4}+z^{\ell }}$
analytic in $\Omega _{1}$, which is an open set containing one of the
branchpoints of $\sqrt{\frac{1}{4}+z^{\ell }}$. This is of course
impossible. So for large $k$, $H_{q_{k}}$ must have a pole in $\Omega _{1}$.%
 
The above argument also works if $\ell $ is odd and we choose 
$$
\Omega _{1}=\Omega _{2}=\{z:r<\left| z\right| <r+\delta \}
$$
for then we cannot find a branch of $\sqrt{\frac{1}{4}+z^{\ell }}$ analytic
in $\Omega _{1}$. (If $\ell $ is even, this argument fails as we may find a
branch analytic \pagebreak in $\Omega _{1}$.)\enddemo

\section{Proof of Theorems 2.4, 2.6}

As a preliminary to proving Theorem 2.4, we rearrange the expression (2.12)
for $\mu _{n}$; recall the notation (3.1).

\nonumproclaim{Lemma 6.1} {\rm (a)} 
\begin{equation}
\mu _{n}(z)=\sum_{j=0}^{[\frac{n+1}{2}]}\frac{q^{-j^{2}}}{(q^{-1};q^{-1})_{j}%
}(zq^{n+1})^{j}\sum_{\ell =0}^{[\frac{n+1}{2}]-j}\frac{z^{\ell }q^{\ell ^{2}}%
}{(q;q)_{\ell }}.
\end{equation}
{\rm (b)} 
\begin{equation}
\nu _{n}(z)=\mu _{n-1}(zq).
\end{equation}
\endproclaim

\demo{Proof} We use the $q$-binomial theorem 
[16, p.~7, eqn.~(1.3.2); p.~236, eq.~(II.4)] 
\begin{equation}
(-uq;q)_{k}=\prod_{j=1}^{k}(1+q^{j}u)=\sum_{j=0}^{k}\left[ 
\begin{array}{c}
k \\ 
j
\end{array}
\right] q^{j(j+1)/2}u^{j}.
\end{equation}
Now,
\begin{eqnarray*}
\left[ 
\begin{array}{c}
n+1-k \\ 
k
\end{array}
\right] &=&\frac{1}{(q;q)_{k}}\prod_{j=1}^{k}(1-q^{j}(q^{n+1-2k}))
\\
&=&\frac{1}{(q;q)_{k}}\sum_{j=0}^{k}\left[ 
\begin{array}{c}
k \\ 
j
\end{array}
\right] q^{j(j+1)/2}(-q^{n+1-2k})^{j}
\\
&=&\sum_{j=0}^{k}\frac{1}{(q;q)_{j}(q;q)_{k-j}}q^{j(j+1)/2+(n+1-2k)j}(-1)^{j}
\end{eqnarray*}
so that from (2.12), 
\begin{eqnarray*}
\mu _{n}(z)&=&\sum_{k=0}^{[\frac{n+1}{2}]}z^{k}q^{k^{2}}\sum_{j=0}^{k}\frac{1}{%
(q;q)_{j}(q;q)_{k-j}}q^{j(j+1)/2+(n+1-2k)j}(-1)^{j}
\\
&
=&\sum_{j=0}^{[\frac{n+1}{2}]}\frac{(-1)^{j}q^{j(j+1)/2+(n+1)j}}{(q;q)_{j}}%
z^{j}\sum_{k=j}^{[\frac{n+1}{2}]}\frac{q^{(k-j)^{2}-j^{2}}}{(q;q)_{k-j}}%
z^{k-j}
\\
&
=&\sum_{j=0}^{[\frac{n+1}{2}]}\frac{q^{-j^{2}}}{(q^{-1};q^{-1})_{j}}%
(zq^{n+1})^{j}\sum_{\ell =0}^{[\frac{n+1}{2}]-j}\frac{q^{\ell ^{2}}}{%
(q;q)_{\ell }}z^{\ell }.
\end{eqnarray*}
In the last line, we used 
\begin{equation}
(q;q)_{j}=(-1)^{j}q^{j(j+1)/2}(q^{-1};q^{-1})_{j}.
\end{equation}

(b) This follows directly from  (2.13). \enddemo

Next, we sketch from \cite[p.~348]{LubinskySaff1987}, \cite[p.~86]
{HardyLittlewood1946} the proof of (2.4):

\nonumproclaim{Lemma 6.2} Let  $q=e^{2\pi i\tau }${\rm ,} with $\tau $  irrational.
The radius of convergence of  $G_{q}$ is 
\begin{equation}
R(q)=\liminf_{n\rightarrow \infty }\left| \prod_{j=1}^{n}(1-q^{j})\right|
^{1/n}=\liminf_{n\rightarrow \infty }| 1-q^{n}|  ^{1/n}.
\end{equation}
\endproclaim

\demo{Proof} Hardy and Littlewood \cite[p.~86]{HardyLittlewood1946}
established the remarkable identity 
\begin{equation}
\sum_{n=0}^{\infty }\frac{z^{n}}{\prod_{j=1}^{n}(1-q^{j})}=\exp
\left(\sum_{n=1}^{\infty }\frac{z^{n}}{n(1-q^{n})}\right)
\end{equation}
and noted that both power series in the last identity have the same radius
of convergence even for our choice of $q$. (The identity was also known to
earlier authors.) The Cauchy-Hadamard formula for the radius of convergence
gives the result.\enddemo

We turn to the

\demo{{P}roof of {\rm (2.15)} and {\rm (2.16)}} Let ${\cal S}$ be an infinite
sequence such that (2.14) holds. We use a section of the sum in the
right-hand side of (6.1). For fixed positive integers $m$, and uniformly for 
$|  z|  \leq r<R(q)$, 
\begin{eqnarray*}
&&\hskip-36pt \lim_{n\rightarrow \infty ,n\in {\cal S}}\sum_{j=0}^{m}\frac{q^{-j^{2}}}{%
(q^{-1};q^{-1})_{j}}(zq^{n+1})^{j}\sum_{\ell =0}^{[\frac{n+1}{2}]-j}\frac{%
q^{\ell ^{2}}}{(q;q)_{\ell }}z^{\ell }\\
&
&\qquad =\sum_{j=0}^{m}\frac{q^{-j^{2}}}{(q^{-1};q^{-1})_{j}}(zq\beta
)^{j}\sum_{\ell =0}^{\infty }\frac{q^{\ell ^{2}}}{(q;q)_{\ell }}z^{\ell
} = G_{q}(z)\overline{\sum_{j=0}^{m}\frac{q^{j^{2}}}{(q;q)_{j}}(\overline{%
zq\beta })^{j}}.
\end{eqnarray*}
Moreover as $m\rightarrow \infty $, the last right-hand side approaches $%
G_{q}(z)\overline{G_{q}(\overline{zq\beta })}$. Moreover, for $|  z| 
\leq r<R(q)$, and uniformly in $n\geq 2m+1$ 
$$
\left| \sum_{j=m+1}^{[\frac{n+1}{2}]}\frac{q^{-j^{2}}}{(q^{-1};q^{-1})_{j}}%
(zq^{n+1})^{j}\sum_{\ell =0}^{[\frac{n+1}{2}]-j}\frac{z^{\ell }q^{\ell ^{2}}%
}{(q;q)_{\ell }}\right|
\leq \left[ \sum_{j=m+1}^{\infty }\frac{r^{j}}{|  (q;q)_{j}|  }\right]
\left[ \sum_{\ell =0}^{\infty }\frac{r^{\ell }}{|  (q;q)_{\ell }|  }%
\right] .
$$
The last right-hand side approaches $0$ as $m\rightarrow \infty $. So we
have (2.15). Then (2.16) follows from the identity (6.2); note that if $%
q^{n}\rightarrow \beta $, then $q^{n-1}\rightarrow \beta /q$.\hfill\qed\enddemo

In the proof of (2.17), we need an identity:

\nonumproclaim{Lemma 6.3} Let  $q=e^{2\pi i\tau }${\rm ,} with $\tau $ irrational.
Then   
\begin{eqnarray}
G_{\overline{q}}(qu)G_{q}(qu)+quG_{\overline{q}}(u)G_{q}(q^{2}u)
&
=&\overline{G_{q}(\overline{qu})}G_{q}(qu)\\
&&+\ qu\overline{G_{q}(\overline{u})}%
G_{q}(q^{2}u)=1.\nonumber
\end{eqnarray}
\endproclaim

\demo{Proof} We first establish the identity 
\begin{equation}
\Gamma :=\sum_{j+k=\ell ;j,k\geq 0}\frac{q^{j^{2}-k^{2}}}{(q;q)_{j}(%
\overline{q};\overline{q})_{k}}u^{j}=q^{\ell ^{2}}u^{\ell }\frac{(q^{1-2\ell
}/u;q)_{\ell }}{(q;q)_{\ell }}.
\end{equation}
Using (6.4), we see that 
\begin{eqnarray*}
\Gamma& =&u^{\ell }\sum_{j+k=\ell ;j,k\geq 0}\frac{q^{j^{2}-k^{2}+k(k+1)/2}}{%
(q;q)_{j}(q;q)_{k}}(-u)^{-k}
\\[4pt]
&=&\frac{q^{\ell ^{2}}u^{\ell }}{(q;q)_{\ell }}\sum_{k=0}^{\ell }\left[ 
\begin{array}{c}
\ell \\[4pt] 
k
\end{array}
\right] q^{k(k+1)/2}(-q^{-2\ell }/u)^{k}
\\[4pt]
&=&q^{\ell ^{2}}u^{\ell }\frac{(q^{1-2\ell }/u;q)_{\ell }}{(q;q)_{\ell }},
\end{eqnarray*}
by the $q$-binomial theorem (6.3). Then by (6.8),  
\begin{eqnarray*}
\Delta &:=&G_{\overline{q}}(qu)G_{q}(qu)+quG_{\overline{q}}(u)G_{q}(q^{2}u)-1
\\[4pt]
&=&\sum_{\ell =1}^{\infty }(uq)^{\ell }\sum_{j+k=\ell ;j,k\geq 0}\frac{%
q^{j^{2}-k^{2}}}{(q;q)_{j}(\overline{q};\overline{q})_{k}}+qu\sum_{\ell
=0}^{\infty }u^{\ell }\sum_{j+k=\ell ;j,k\geq 0}\frac{q^{j^{2}-k^{2}}q^{2j}}{%
(q;q)_{j}(\overline{q};\overline{q})_{k}}
\\[4pt]
&=&\sum_{\ell =1}^{\infty }(uq)^{\ell }q^{\ell ^{2}}\frac{(q^{1-2\ell
};q)_{\ell }}{(q;q)_{\ell }}+qu\sum_{\ell =0}^{\infty }u^{\ell }q^{\ell
^{2}}q^{2\ell }\frac{(q^{1-2\ell -2};q)_{\ell }}{(q;q)_{\ell }}
\\[4pt]
&=&\sum_{\ell =1}^{\infty }u^{\ell }q^{\ell ^{2}}\frac{(q^{1-2\ell };q)_{\ell
-1}}{(q;q)_{\ell -1}}\left[ q^{\ell }\frac{1-q^{-\ell }}{1-q^{\ell }}%
+1\right] =0.
\\
\noalign{\vskip-36pt}
\end{eqnarray*}
\enddemo

\phantom{howIwish}

We turn to the

\demo{Proof of {\rm (2.17)}} We assume that ${\cal S}$ satisfies (2.14). From
(3.6), 
$$
\frac{H_{q}(z)-\frac{\mu _{n}(z)}{\nu _{n}(z)}}{%
(-1)^{n}z^{n+1}q^{(n+1)(n+2)/2}}=\frac{1}{\nu _{n}(z)[\nu
_{n}(z)H_{q}(q^{n+1}z)+q^{n+1}z\nu _{n-1}(z)]}.
$$
Here $H_{q}(q^{n+1}z)\rightarrow H_{q}(q\beta z)$ as $n\rightarrow \infty $
through ${\cal S}$, uniformly in compact subsets of $\left| z\right|
<R\left( q\right) $ omitting zeros of $G_{q}(q^{2}\beta z)$. Using (2.16),
we obtain uniformly in compact subsets of $|  z|  <R(q)$ omitting zeros
of $G_{q}(q^{2}\beta z)$, 
\begin{eqnarray*}&&
\lim_{n\rightarrow \infty ,n\in {\cal S}}\nu _{n}(z)\left[ \nu
_{n}(z)H_{q}(q^{n+1}z)+q^{n+1}z\nu _{n-1}(z)\right]
\\[4pt] &&\qquad
=G_{q}(qz)\overline{G_{q}(\overline{\beta qz})}\left[ G_{q}(qz)\overline{%
G_{q}(\overline{\beta qz})}H_{q}(q\beta z)+q\beta zG_{q}(qz)\overline{G_{q}(%
\overline{\beta z})}\right]
\end{eqnarray*}\pagebreak
\begin{eqnarray*} &&\qquad
=\frac{G_{q}(qz)^{2}\overline{G_{q}(\overline{\beta qz})}}{G_{q}(q^{2}\beta
z)}\left[ \overline{G_{q}(\overline{\beta qz})}G_{q}(\beta qz)+q\beta z%
\overline{G_{q}(\overline{\beta z})}G_{q}(q^{2}\beta z)\right] \\[4pt] &&\qquad=\frac{%
G_{q}(qz)^{2}\overline{G_{q}(\overline{\beta qz})}}{G_{q}(q^{2}\beta z)}
\end{eqnarray*}
by Lemma 6.3. The result then follows.\enddemo

Now we can turn to the\newline
\demo{Proof of Corollary {\rm 2.5}} Note first that $\mu _{n},\nu _{n}$ have no
common zeros (see (3.4)), so every zero of $\nu _{n}$ is a pole of $\mu
_{n}/\nu _{n}$. But the limit relation (2.16) shows that for the appropriate
subsequence, if $H_{q}$ has $\ell $ poles on $|  z|  =r$, that is $%
G_{q}(qz)$ has $\ell $ zeros there, then for large $n\in S$, $\nu _{n}$ has
at least $2\ell $ zeros counting multiplicity in any neighbourhood of that
circle. As every subsequence of positive integers contains a subsequence
satisfying (2.14), for some $\beta $, the result follows.\enddemo

Our proof of Theorem 2.6 is very similar to that of Theorem 2.4 in \cite{LubinskySaff1987}, but we provide most of the
details.

\demo{Proof of Theorem {\rm 2.6}} Let us set 
$$
H_{j}(u):=\prod_{k=1}^{j}(1+q^{-k}u).
$$
We see that 
$$
\left[ 
\begin{array}{c}
n \\ 
j
\end{array}
\right] =H_{j}(-q^{n+1})/\prod_{k=1}^{j}(1-q^{k}).
$$
Since $\{-q^{n+1}\}_{n=j}^{\infty }$ is dense on the unit circle, we see
that 
$$
\Gamma _{j}:=\sup_{n\geq j}\left| \left[ 
\begin{array}{c}
n \\ 
j
\end{array}
\right] \right| =  \Vert H_{j}\Vert _{L_{\infty }(|  z|  =1)}/\left| 
\prod_{k=1}^{j}(1-q^{k}) \right| .
$$
Then we see that 
\begin{equation}
\Vert \mu _{n}\Vert _{L_{\infty }(|  z|  \leq r)}\ \leq\
\sum_{j=0}^{\infty }\Gamma _{j}r^{j},
\end{equation}
with a similar inequality for $\nu _{n}$. It is shown in \cite[p.\ 345]
{LubinskySaff1987} with the aid of the theory of uniform distribution (there 
$q$ is replaced by $q^{-1}$) that 
$$
\lim_{j\rightarrow \infty }\Vert H_{j}\Vert _{L_{\infty }(|  z| 
=1)}^{1/j}=1
$$
so that
$$
\limsup_{j\rightarrow \infty }\Gamma _{j}^{1/j}=1/R(q).
$$
Thus the series in (6.9) converges if $r<R(q)$. Since the series is
independent of $n$, we have the required uniform boundedness of $\{\mu
_{n}\},\{\nu _{n}\}$. Next if $r$ is such that 
$$
C:=\sup_{n\geq 1}\Vert \nu _{n}\Vert _{L_{\infty }(|  z|  \leq
r)}<\infty ,
$$
Cauchy's inequalities give for $\left[ \frac{n}{2}\right] \geq k,$%
$$
\left|  \left[ 
\begin{array}{c}
n-k \\ 
k
\end{array}
\right] \right| \leq\ \Vert \nu _{n}\Vert _{L_{\infty }(|  z|  \leq
r)}/r^{k}\leq C/r^{k}.
$$
As $n-k$ may assume any integral value $\geq k$ as $n$ runs over integers with 
$\left[ \frac{n}{2}\right] \geq k$, we obtain 
$$
\Gamma _{k}\leq C/r^{k}.
$$
Taking $k^{\rm th}$ roots and letting $k\rightarrow \infty $ give 
$$
1/R(q)\leq 1/r
$$
so that  $r\leq R(q)$. Thus $\{\mu _{n}\},\{\nu _{n}\}$ cannot be uniformly
bounded in $|  z|  \leq r$ if $r>R(q).$
\enddemo

\section{Proof of Theorem 1.1}

We shall assume throughout that we are dealing with $q$ on the unit circle
such that $R\left( q\right) =1$. (Later on in this section, we shall
specialize $q$ to that given in Theorem 1.1.) Note that by Theorem 2.2, $%
H_{q}$ is meromorphic in the unit ball, and $G_{q}$ is analytic in the unit
ball, and both have a natural boundary on the unit circle. We first outline
the main steps in the proof of Theorem 1.1:
\vglue8pt
(I) We again show how a counterexample to the Baker-Gammel-Wills Conjecture
follows if  for all $\beta $ on the unit circle, 
\begin{equation}
\left\{ z:G_{q}\left( qz\right) =0\right\} \neq \left\{ z:G_{q}\left( 
\overline{\beta qz}\right) =0\right\} .
\end{equation}
\vglue8pt (II) We show that to prove (7.1), instead of considering the full series for 
$G_{q}$, it suffices to show something like (7.1) for a partial sum 
$$
S_{m,q}\left( z\right) :=\sum_{j=0}^{m}\frac{q^{j^{2}}}{\left( q;q\right)
_{j}}z^{j}.
$$
(We shall use $m=50$.) This is achieved using Rouche's theorem and an
estimate for the tail $G_{q}-S_{50,q}$.\vglue8pt

(III) We use the principal of the argument to count the number of zeros of $%
S_{50,q}$ inside certain circles. To evaluate the integrals counting the
zeros, we use an elementary integration rule, and establish a rigorous
estimate for the error. The calculation is performed using Matlab 6.0.
Moreover, we use Matlab 6.0 and Mathematica 3.0 to estimate below the
minimum modulus of $S_{50,q}$ on certain circles. Our numerical evaluations
involve only evaluating $S_{50,q}$ and its first three derivatives at
various explicit points, and sums involving it, or comparisons of its
absolute value at a definite set of points. The calculation does not make
any use of ``black-box'' zero finding routines, or numerical quadrature
routines. These calculations show that the two zeros of $S_{50,q}$ of
smallest modulus have the desired asymmetry, and we then deduce the same for 
$G_{q}$.\vglue8pt

We now turn to these steps in detail:

\demo{{\rm (I)} It suffices to prove that {\rm (7.1)} holds for every $\left| \beta
\right| =1$}
Let us suppose that (7.1) holds for every $\left| \beta \right| =1$, but
that some subsequence $\left\{ \mu _{n}/\nu _{n}\right\} _{n\in {\cal S}}$
has 
$$
\lim_{n\rightarrow \infty ,n\in {\cal S}}\mu _{n}/\nu _{n}=H_{q}
$$
uniformly in compact subsets of the unit ball omitting poles of $H_{q}$. By
extracting a further subsequence, we may assume that for some $\left| \beta
\right| =1,$%
$$
\lim_{n\rightarrow \infty ,n\in {\cal S}}q^{n}=\beta .
$$
Now let $K$ be a closed ball in the unit ball, containing in its interior 
zeros of $G_{q}\left( \overline{\beta qz}\right) $ that are not zeros of $%
G_{q}\left( qz\right) $. We may also assume that $K$ does not contain any
zeros of $G_{q}\left( qz\right) $, that is, poles of $H_{q}$. For large $%
n\in {\cal S}$, $\mu _{n}/\nu _{n}$ has no poles in $K$, because of the
assumed uniform convergence. Since $\mu _{n}$ and $\nu _{n}$ are coprime
polynomials (recall (3.4)), this forces $\nu _{n}$ not to have zeros in $K$
for large $n$. But that contradicts the uniform convergence in (2.16), which
by Hurwitz' Theorem, shows that each zero of $G_{q}\left( \overline{\beta qz}%
\right) $ in $K$ must attract zeros of $\nu _{n},n\in {\cal S}$.
\enddemo

Rather than working with (7.1), it will be easier to work with:

\demo{Reformulation of {\rm (7.1):} The zeros of $G_{q}$  are not
symmetric about any line through $0$}
That is, after reflecting the set of zeros about any line through $0$, we
obtain a different set.

To see this, observe that if, for example, $\beta =1$, (7.1) requires that
the zeros of $G_{q}$ do not occur in conjugate pairs. That is, the zeros of $%
G_{q}$ are not symmetric about the real line. Of course the case of general $%
\beta $ is similar.
\enddemo
 {\rm (II)} {\it Estimation of the tail $G_{q}-S_{50,q}$}.
We estimate the tail for a special class of $q$ including that in Theorem
1.1.

\nonumproclaim{Lemma 7.1} Let  $q=e^{2\pi i\tau }${\rm ,} where for some positive integer $%
c\geq 2${\rm ,}
\begin{equation}
\tau :=\frac{1|}{|c}+\frac{1|}{|1}+\frac{1|}{|1}+\frac{1|}{|1}+\cdots =\frac{1}{%
c+\frac{1}{2}\left( \sqrt{5}-1\right) }.
\end{equation}
Then for $0<r<s<1,$%
\begin{eqnarray}
&&\max_{\left| z\right| \leq r}\left| G_{q}\left( z\right) -S_{m,q}\left(
z\right) \right| \\
&&\qquad\quad\leq \frac{\left( r/s\right) ^{m+1}}{1-r/s}\exp \left( S_{0}+\frac{\sqrt{5}%
}{4\left( 1-\alpha ^{-8}\right) }\frac{s^{2c+1}}{1-s}\right) =:T, \nonumber 
\end{eqnarray}
where 
\begin{equation}
S_{0}:=\sum_{n=1}^{2c}\frac{s^{n}}{2n\left| \sin n\pi \tau \right| }%
;\;\alpha :=\frac{1}{2}\left( \sqrt{5}+1\right) .
\end{equation}
\endproclaim

\demo{Proof} 
From (6.6), and Cauchy's estimates, for $0<s<1,$%
\begin{eqnarray}
\frac{1}{\left| \left( q;q\right) _{n}\right| } &\leq &s^{-n}\max_{\left|
z\right| =s}\left| \exp \left( \sum_{n=1}^{\infty }\frac{z^{n}}{n\left(
1-q^{n}\right) }\right) \right| \\
&\leq &s^{-n}\exp \left( \sum_{n=1}^{\infty }\frac{s^{n}}{2n\left| \sin n\pi
\tau \right| }\right) . \nonumber 
\end{eqnarray}

To estimate the last series, we use the classical continued fraction
expansion (7.2) of $\tau $. Let $\left\{ \pi _{j}/\chi _{j}\right\}
_{j=0}^{\infty }$ be the convergents in the continued fraction (7.2) of $%
\tau $. With initial conditions 
\begin{eqnarray*}
\pi _{0} &=&0;\qquad \pi _{1}=1; \\
\chi _{0} &=&1;\qquad \chi _{1}=c,
\end{eqnarray*}
they satisfy the following recurrence relations for $n\geq 2$ (\cite
{HardyWright1975}, \cite{Lang1966}): 
\begin{eqnarray}
\pi _{n} &=&\pi _{n-1}+\pi _{n-2};  \\
\chi _{n} &=&\chi _{n-1}+\chi _{n-2}. \nonumber
\end{eqnarray}
Of course this special form of the recurrence relation exists because all
partial quotients in (7.2), other than the first, are unity. To solve this
constant coefficient difference equation, one uses the characteristic
equation 
$$
x^{2}-x-1=0,
$$
with roots 
$$
\alpha :=\frac{1}{2}\left( 1+\sqrt{5}\right) ;\qquad \beta :=\frac{1}{2}\left( 1-%
\sqrt{5}\right) .
$$
Standard methods and the initial conditions give for $n\geq 0,$%
\begin{eqnarray}
\pi _{n} &=&\left( \alpha ^{n}-\beta ^{n}\right) /\sqrt{5}; \\
\chi _{n} &=&\left( \left( c-\beta \right) \alpha ^{n}+\left( \alpha
-c\right) \beta ^{n}\right) /\sqrt{5}.
\end{eqnarray}
Moreover, the form of $\tau $ and a simple calculation give 
\begin{equation}
\chi _{n}\tau -\pi _{n}=\frac{\beta ^{n}}{c-\beta }=\frac{\left( -\alpha
\right) ^{-n}}{c-\beta },
\end{equation}
so that 
\begin{equation}
\Vert \chi _{n}\tau \Vert =\frac{1}{\left( c-\beta \right) \alpha
^{n}}.
\end{equation}
Note that the last right-hand side is bounded above by 
$$
\frac{1}{\left( 1-\beta \right) \alpha }=\frac{1}{\alpha ^{2}}<\frac{1}{2},
$$
so that  (7.10) is true even for $n=1$.

Next, we fix $j\geq 3$, and consider $n$ such that 
$$
\chi _{j}\leq n<\chi _{j+1}.
$$
By the best approximation property of continued fractions \cite{Lang1966}, 
$$
\Vert n\tau \parallel\ \geq\ \Vert \chi _{j}\tau \Vert $$
and hence 
\begin{equation}
n\Vert n\tau \parallel\ \geq \chi _{j}\Vert \chi _{j}\tau \parallel\ =%
\frac{1}{\sqrt{5}}\left( 1+\frac{\alpha -c}{c-\beta }\frac{\left( -1\right)
^{j}}{\alpha ^{2j}}\right) ,
\end{equation}
by (7.8) and (7.10). Here 
\begin{equation}
0<\left( -1\right) \frac{\alpha -c}{c-\beta }=\frac{c-\alpha }{c+\alpha ^{-1}%
}<1.
\end{equation}
Thus if $j=3$, the right-hand side of (7.11) exceeds $1/\sqrt{5}$. If $j\geq
4$, it exceeds 
$$
\frac{1}{\sqrt{5}}\left( 1-\alpha ^{-8}\right) .
$$
Then for $n\geq \chi _{3}$, 
\begin{eqnarray*}
\frac{1}{2n\left| \sin n\pi \tau \right| } &=&\frac{1}{2n\left( \sin \pi
\Vert n\tau \Vert \right) } \\
&\leq &\frac{1}{4n\Vert n\tau \Vert }\leq \frac{\sqrt{5}}{4\left(
1-\alpha ^{-8}\right) }.
\end{eqnarray*}
Since $\chi _{3}=2c+1$, we deduce from (7.5) that 
$$
\frac{1}{\left| \left( q;q\right) _{n}\right| }\leq s^{-n}\exp \left( S_{0}+%
\frac{\sqrt{5}}{4\left( 1-\alpha ^{-8}\right) }\frac{s^{2c+1}}{1-s}\right) ,
$$
where $S_{0}$ is given by (7.4). Multiplying this by $r^{n}$ and adding over 
$n\geq m+1$ gives (7.3). \enddemo

We shall also need an immediate consequence of Rouche's theorem and Lemma
7.1:

\nonumproclaim{Lemma 7.2} Assume the hypotheses and notation of Lemma {\rm 7.1.} Let  $\gamma $%
  be a simple closed curve in  $\left\{ z:\left| z\right| \leq
r\right\} $  such that 
\begin{equation}
\min_{z\in \gamma }\left| S_{m,q}\left( z\right) \right| >T.
\end{equation}
 Then $S_{m,q}$  and  $G_{q}$  have the same total
multiplicity of zeros inside  $\gamma $.
\endproclaim

 (III) {\it Use of Mathematica {\rm 3.0} and Matlab {\rm 6} to verify what is needed
for }$S_{50,q}$.
In applying Lemma 7.2, we need to estimate the minimum modulus of $S_{m,q}$.
Later on, we shall also need to estimate the maximum modulus of $%
S_{50,q}^{(j)},j=0,1,2,3$. Since we prefer to evaluate $S_{m,q}$ only at a
definite set of points, rather than relying on a ``black-box'' to find a
minimum, we need:
\nonumproclaim{Lemma 7.3} Let $P$ be a polynomial of degree  $\leq n${\rm ,} let  $%
m\geq 1${\rm ,} and let $\gamma $  be the circle $\left\{
z:\left| z-a\right| =\varepsilon \right\} $. Then 
\begin{eqnarray}
\min_{z\in \gamma }\left| P\left( z\right) \right|  
&\geq &\min_{1\leq j\leq m}\left| P\left( a+\varepsilon e^{2\pi ij/m}\right)
\right|\\[4pt]
&& -\frac{\pi n}{m}\max_{1\leq j\leq 2n}\left| P\left( a+\varepsilon
e^{2\pi ij/\left( 2n\right) }\right) \right| ;\nonumber
\\[4pt]
 \max_{z\in \gamma }\left| P\left( z\right) \right|  
&\leq &\max_{1\leq j\leq m}\left| P\left( a+\varepsilon e^{2\pi ij/m}\right)
\right| \\[4pt]
&&+\frac{\pi n}{m}\max_{1\leq j\leq 2n}\left| P\left( a+\varepsilon
e^{2\pi ij/\left( 2n\right) }\right) \right| .\nonumber
\end{eqnarray}
\endproclaim
 
\demo{Proof} 
We may assume that $a=0$ and $\varepsilon =1$, so that $\gamma $ is the unit
circle. The general case follows by replacing $P\left( z\right) $ by $%
P\left( a+\varepsilon z\right) $. We use a Duffin-Schaefer type inequality
due to Frappier, Rahman, and Ruscheweyh \cite{Frappieretal1985}, \cite[p.~690]
{Milovanovicetal1994}: 
\begin{equation}
\Vert P^{\prime }\Vert _{L_{\infty }\left( \left| z\right| =1\right)
}\leq n\max_{1\leq j\leq 2n}\left| P\left( e^{2\pi ij/\left( 2n\right)
}\right) \right| .
\end{equation}
Let $z=e^{i\theta }$, $\theta \in \left[ 0,2\pi \right] $, and choose $j\in
\left\{ 0,1,2,\ldots,m\right\} $ such that $\left| \theta -\frac{2\pi j}{m}%
\right| $ is as small as is positive. Then 
\begin{eqnarray*}
\left| P\left( z\right) -P\left( e^{2\pi ij/m}\right) \right| &=&\left|
\int_{2\pi j/m}^{\theta }P^{\prime }\left( e^{it}\right) ie^{it}dt\right| \\
&\leq &\left| \theta -\frac{2\pi j}{m}\right| \Vert P^{\prime }\Vert _{L_{\infty }\left( \left| z\right| =1\right) } \\
&\leq &\frac{\pi }{m}\left( n\max_{1\leq k\leq 2n}\left| P\left(
 e^{2\pi ik/\left( 2n\right) }\right) \right| \right) ,
\end{eqnarray*}
by (7.16). Then (7.14) and (7.15) follow.\enddemo

To count the number of zeros of $S_{50,q}$, we use the principal of the
argument: if $\gamma $ is a circle $\left\{ z:\left| z-a\right| =r\right\} $
on which $S_{50,q}$ has no zeros, 
$$
I\left( \gamma \right) :=\frac{1}{2\pi i}\int_{\gamma }\frac{%
S_{50,q}^{\prime }}{S_{50,q}}=\frac{r}{2\pi }\int_{-\pi }^{\pi }e^{i\theta }%
\frac{S_{50,q}^{\prime }}{S_{50,q}}\left( a+re^{i\theta }\right) \ d\theta
$$
is the total multiplicity of zeros of $S_{50,q}$ inside $\gamma $. We
approximate $I\left( \gamma \right) $ by the simple rule 
$$
I_{m}\left( \gamma \right) =\frac{r}{m}\sum_{j=0}^{m-1}e^{2\pi ij/m}\frac{%
S_{50,q}^{\prime }}{S_{50,q}}\left( a+re^{2\pi ij/m}\right) .
$$
Following is an elementary estimate for the error. It is by no means
optimal, but suffices for our purposes:
\nonumproclaim{Lemma 7.4} 
\begin{equation}
\left| I\left( \gamma \right) -I_{m}\left( \gamma \right) \right| \leq \frac{%
r\pi \sqrt{2}}{m^{2}}\Phi ,
\end{equation}
 where 
\begin{eqnarray}&&\\
\Phi &:=&\frac{\max_{\gamma }\left| S_{50,q}^{\prime }\right| }{\min_{\gamma
}\left| S_{50,q}\right| }+3r\frac{\max_{\gamma }\left| S_{50,q}^{\prime
\prime }\right| }{\min_{\gamma }\left| S_{50,q}\right| }+3r\left( \frac{%
\max_{\gamma }\left| S_{50,q}^{\prime }\right| }{\min_{\gamma }\left|
S_{50,q}\right| }\right) ^{2}  +r^{2}\frac{\max_{\gamma }\left| S_{50,q}^{\prime \prime \prime }\right| }{%
\min_{\gamma }\left| S_{50,q}\right| }\nonumber \\
&&+\ 3r^{2}\frac{\max_{\gamma }\left|
S_{50,q}^{\prime }\right| \max_{\gamma }\left| S_{50,q}^{\prime \prime
}\right| }{\left( \min_{\gamma }\left| S_{50,q}\right| \right) ^{2}}%
+2r^{2}\left( \frac{\max_{\gamma }\left| S_{50,q}^{\prime }\right| }{%
\min_{\gamma }\left| S_{50,q}\right| }\right) ^{3}.\nonumber
\end{eqnarray}
\endproclaim

\demo{Proof} 
First recall that if 
$$
R\left( z\right) =\sum_{j=-\left( m-1\right) }^{m-1}c_{j}z^{j}
$$
is a trigonometric polynomial of degree $<m$, then 
$$
\frac{1}{2\pi }\int_{-\pi }^{\pi }R\left( e^{i\theta }\right) d\theta =\frac{%
1}{m}\sum_{j=0}^{m-1}R\left( e^{2\pi ij/m}\right) .
$$
Next, for any continuous complex-valued function $g$ defined on the unit
circle, we deduce that 
\begin{eqnarray*}
E &:=&\left| \frac{1}{2\pi }\int_{-\pi }^{\pi }g\left( e^{i\theta }\right)
d\theta -\frac{1}{m}\sum_{j=0}^{m-1}g\left( e^{2\pi ij/m}\right) \right| \\
&\leq &2\inf \Vert g-R\Vert _{L_{\infty }\left( \left| z\right|
=1\right) },
\end{eqnarray*}
where the inf is taken over all trigonometric polynomials $R$ of degree $<m$%
. To estimate this error of approximation, we use one case of the Favard
estimates \cite[Thm.\ 4.3, p.\ 214]{DeVoreLorentz1993}. This asserts that if $%
h $ is a twice continuously differentiable real valued function defined on $%
\left[ -\pi ,\pi \right] $ , there exists a trigonometric polynomial $R$ of
degree $<m$ such that 
$$
\Vert h-R\Vert _{L_{\infty }\left[ -\pi ,\pi \right] }\ \leq \frac{\pi 
}{2m^{2}}\Vert h^{\prime \prime }\Vert _{L_{\infty }[-\pi ,\pi ]}.
$$
Applying this to $h\left( \theta \right) ={\rm Re}g\left( e^{i\theta
}\right) $ and $h\left( \theta \right) ={\rm Im}g\left( e^{i\theta }\right) 
$ gives, if the second derivatives of $g$ are continuous, 
$$
E\leq \frac{\sqrt{2}\pi }{m^{2}}\sup \left\{ \left| \frac{d^{2}}{d\theta ^{2}%
}g\left( e^{i\theta }\right) \right| :\theta \in \left[ -\pi ,\pi \right]
\right\} .
$$
(The $\sqrt{2}$ can probably be dropped, as the Favard estimate probably
also applies to complex-valued functions.) Applying this to 
$$
g\left( e^{i\theta }\right) :=re^{i\theta }\frac{S_{50,q}^{\prime }}{S_{50,q}%
}\left( a+re^{i\theta }\right) ,
$$
we obtain 
$$
\left| I\left( \gamma \right) -I_{m}\left( \gamma \right) \right| \leq \frac{%
r\sqrt{2}\pi }{m^{2}}\sup \left\{ \left| \frac{d^{2}}{d\theta ^{2}}\left[
e^{i\theta }\frac{S_{50,q}^{\prime }}{S_{50,q}}\left( a+re^{i\theta }\right)
\right] \right| :\theta \in \left[ -\pi ,\pi \right] \right\} .
$$
An explicit calculation of this derivative in terms of the derivatives of $%
S_{50,q}$ and some elementary estimates then yield (7.18). \enddemo

Now let us choose in Lemma 7.1,\begin{itemize}
\item[(i)] $c=m=50$, so that $\tau $ of (7.2) becomes $\tau $ of Theorem 1.1.%
\item[(ii)] $s=0.9$ and $r=0.46.$\end{itemize}
A Mathematica 3.0 calculation shows that $T$ of (7.3) is given by 
$$
T=1.04093\cdots\times 10^{-10}.
$$
Thus 
\begin{equation}
\max_{\left| z\right| \leq 0.46}\left| G_{q}-S_{50,q}\right| \left( z\right)
\leq T<1.04094\times 10^{-10}.
\end{equation}
Next, let 
\begin{eqnarray*}
z_{1} &:=&-0.299076+0.145052i; \\
z_{2} &:=&-0.269527+0.306036i
\end{eqnarray*}
so that 
\begin{eqnarray*}
\left| z_{1}\right| &=&0.332395\cdots; \\
\left| z_{2}\right| &=&0.407802\cdots\;.
\end{eqnarray*}
Also, let 
\begin{eqnarray*}
\gamma _{0} &:=&\left\{ t:\left| t\right| =0.46\right\}, \\
\gamma _{j} &:=&\left\{ t:\left| t-z_{j}\right| =0.01\right\} ,\qquad j=1,2.
\end{eqnarray*}
Mathematica 3.0 calculated $z_{1}$ and $z_{2}$ as zeros of $S_{50,q}$, but
we do not need in our proof to assume that these are (approximations to)
zeros of $S_{50,q}$.

Now let us summarize what we need. Suppose that we can show \vglue4pt
(A) 
\begin{equation}
\min_{\gamma _{j}}\left| S_{50,q}\right| >T,\qquad j=0,1,2.
\end{equation}
\vglue4pt
(B) For appropriate choices of $m_{j},j=0,1,2,$%
\begin{equation}
\left| I\left( \gamma _{j}\right) -I_{m_{j}}\left( \gamma _{j}\right)
\right| <\frac{1}{2},\qquad  j=0,1,2.
\end{equation}
\vglue4pt
(C) The closest integer to $I_{m_{j}}\left( \gamma _{j}\right) $ is $2$ for $%
j=0$ and $1$ for $j=1,2$.

Then (A) and (7.19) show that 
$$
\min_{\gamma _{j}}\left| S_{50,q}\right| >\max_{\gamma _{j}}\left|
G_{q}-S_{50,q}\right| ,
$$
so that by Lemma 7.2, $S_{50,q}$ and $G_{q}$ have the same number of zeros
inside $\gamma _{j},j=0,1,2$. Next, (B) shows that the closest integer to $%
I_{m_{j}}\left( \gamma _{j}\right) $ is $I\left( \gamma _{j}\right) ,\qquad j=0,1,2$%
. Then (C) shows that 
$$
I\left( \gamma _{0}\right) =2,
$$
and 
$$
I\left( \gamma _{1}\right) =1=I\left( \gamma _{2}\right) .
$$
Then $S_{50,q}$ and hence $G_{q}$ have zeros of total multiplicity $2$
inside $\gamma _{0}$, and these must be the simple zeros inside $\gamma
_{j},j=1,2$. If we denote these zeros by $u_{j},j=1,2$, then 
$$
\left| u_{j}-z_{j}\right| <0.01,\qquad j=1,2,
$$
so  that 
\begin{eqnarray*}
0 &<&\left| u_{1}\right| <\left| u_{2}\right| ; \\
\arg \left( u_{1}\right) &\neq &\arg \left( u_{2}\right) .
\end{eqnarray*}
The figure below contains the curves $\gamma _{j},j=0,1,2$. The asymmetry of
the zeros $z_{1},z_{2}$ and the circles containing the corresponding zeros $%
u_{1},u_{2}$ of $G_{q}$ is then clear. Thus the zeros of $G_{q}$ cannot be
symmetric about any line through $0$, and in fact, this is true even for the
two zeros $u_{1}$ and $u_{2}$ of smallest modulus. So we have satisfied our
reformulation of (7.1) and $H_{q}$ serves as a counterexample to the
Baker-Gammel-Wills Conjecture.

Even more, it follows by the exact same argument from Step I above, that
since the zeros $u_{1}$ and $u_{2}$ of smallest modulus cannot be symmetric
about any line through $0$, we cannot have a subsequence of $\left\{ \mu
_{n}/\nu _{n}\right\} _{n=1}^{\infty }$ converging uniformly in all compact
subsets of $\left\{ z:\left| z\right| <0.46\right\} $. So we have completed
the proof of Theorem 1.1. 
\figin{fig1}{600}
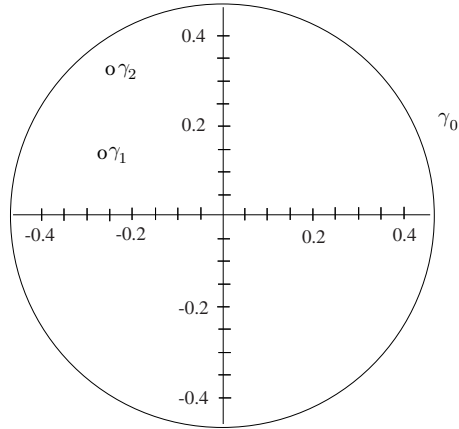
\centerline{Smallest zeros of $G_q$}

\vglue12pt

Now we turn to verifying (A), (B), (C). Initially calculations were
performed using Mathematica 3.0, and later using Matlab 4. The calculations
for the final version of this manuscript were performed using Matlab 6.

\demo{{P}roof of {\rm (A)} for $j=0$}
A Matlab 6 calculation shows that 
\begin{eqnarray*}
\max_{1\leq j\leq 10^{5}}\left| S_{50,q}\left( 0.46e^{2\pi ij/10^{5}}\right)
\right| &=&37.9643\cdots; \\
\min_{1\leq j\leq 2\times 10^{7}}\left| S_{50,q}\left( 0.46e^{2\pi ij/\left(
2\times 10^{7}\right) }\right) \right| &=&0.01307\cdots\;.
\end{eqnarray*}
Then Lemma 7.3 with $m=2\times 10^{7}$ and $n=50$ gives 
$$
\min_{\left| z\right| =0.46}\left| S_{50,q}\left( z\right) \right| \geq
0.01277\cdots>T.
$$
(Of course the maximum above exceeds that needed for Lemma 7.3).\enddemo

{\it Proof of} (A) {\it for }$j=1$.
A Matlab 6 calculation shows that 
\begin{eqnarray*}
\max_{1\leq j\leq 10^{5}}\left| S_{50,q}\left( z_{1}+0.01e^{2\pi
ij/10^{5}}\right) \right| &=&0.04735\cdots; \\
\min_{1\leq j\leq 10^{5}}\left| S_{50,q}\left( z_{1}+0.01e^{2\pi
ij/10^5}\right) \right| &=&0.03581\dots\;.
\end{eqnarray*}
Then Lemma 7.3 with $m=10^{5}$ and $n=50$ gives 
$$
\min_{\left| z-z_{1}\right| =0.01}\left| S_{50,q}\left( z\right) \right|
\geq 0.03580\cdots>T.
$$
\demo{{P}roof of {\rm (A)} for $j=2$}
A Matlab 6 calculation shows that 
\begin{eqnarray*}
\max_{1\leq j\leq 10^{5}}\left| S_{50,q}\left( z_{2}+0.01e^{2\pi
ij/10^{5}}\right) \right| &=&0.01516\cdots; \\
\min_{1\leq j\leq 10^{5}}\left| S_{50,q}\left( z_{2}+0.01e^{2\pi
ij/10^5}\right) \right| &=&0.01205\cdots\;.
\end{eqnarray*}
Then Lemma 7.3 with $m=10^{5}$ and $n=50$ gives 
$$
\min_{\left| z-z_{2}\right| =0.01}\left| S_{50,q}\left( z\right) \right|
\geq 0.01203\cdots>T.
$$
\newline
Thus we have completed the proof of (A). Next, we turn to (B) and (C). We
use Lemma 7.4, and so must estimate $\Phi $ from (7.18). To estimate $%
\max_{\gamma _{j}}\left| S_{50,q}^{(\ell )}\right| $, $\ell =0,1,2,3$, we
used (7.15) of Lemma 7.3 with $n=50^{\phantom{|}}\hskip-4pt$ and $m=10^{5}$ in all cases, applied to 
$S_{50,q}^{(\ell )}$. The results appear in the table below: 
$$
\begin{tabular}{|l|l|l|l|l|l|}
\hline
$j$ & $\min_{\gamma _{j}}\left| S_{50,q}\right| $ & $\max_{\gamma
_{j}}\left| S_{50,q}\right| $ & $\max_{\gamma _{j}}\left| S_{50,q}^{\prime
}\right| $ & $\max_{\gamma _{j}}\left| S_{50,q}^{\prime \prime }\right| $ & $%
\max_{\gamma _{j}}\left| S_{50,q}^{\prime \prime \prime }\right| $ \\ \hline
$0$ & $0.01277..$ & $38.021..$ & $311.06..$ & $2,562.2..$ & $21,278.9..$ \\ 
\hline
$1$ & $0.03580..$ & $0.04742..$ & $5.3911..$ & $137.77..$ & $2408.95..$ \\ 
\hline
$2$ & $0.01203..$ & $0.01519..$ & $1.7003..$ & $39.4139..$ & $941.25..$ \\ 
\hline
\end{tabular}
$$
If we substitute these estimates in Lemma 7.4 with $m_{j}=2\times 10^{7}$
for $j=0$, and $m_{j}=10^{5}$ for $j=1,2,$ we obtain the estimates in the
third column of the following table. 
$$
\begin{tabular}{|l|l|l|l|}
\hline
$j$ & $m_{j}$ & $\left| I\left( \gamma _{j}\right) -I_{m_{j}}\left( \gamma
_{j}\right) \right| \leq $ & $I_{m}\left( \gamma _{j}\right) $ \\ \hline
$0$ & $2\times 10^{7}$ & $3.17..\times 10^{-2}$ & $2-\left( 1.5\times
10^{-13}+10^{-14}i\right) $ \\ \hline
$1$ & $10^{5}$ & $8.03..\times 10^{-9}$ & $1+10^{-14}$ \\ \hline
$2$ & $10^{5}$ & $8.64..\times 10^{-9}$ & $1$ \\ \hline
\end{tabular}
$$
We see that the error in numerical integration is far less than $\frac{1}{2}$%
. Furthermore, to at least 12 decimals, $I_{m_{j}}\left( \gamma _{j}\right) $
is $2$ for $j=0$, and $1$ for $j=1,2$. So we have completed (B), (C) and the
proof of Theorem 1.1.\hfill\qed\enddemo

{\it Remarks}.
(i) Obviously the accuracy of the calculation is fundamental in the above
proof. We tested our calculation of $S_{50,q}$ in Mathematica 3.0 by
calculating it in two different ways. The second was more careful, involving
separating out real and imaginary, and then positive and negative parts of
the summand. Adding positive parts separately and then negative parts
separately to avoid roundoff, we obtained agreement with the simpler method
of calculation to 12 decimals. Similar checking of other calculations
indicated accuracy to 12 decimals, and independent calculations on Matlab
gave the same results.\pagebreak

(ii) We again emphasize that the calculations above involved only evaluation
of $S_{50,q}$ and its first three derivatives at translated roots of unity,
together with their absolute values, and sums involving them. No ``black-box'' routines were involved. In the first version of this paper, the author
used the ``black-box'' routines of Mathematica and Matlab to evaluate $%
I\left( \gamma _{j}\right) $, and did not provide an error estimate. The
author must thank the referee for requesting an error estimate, which led to
the inclusion of Lemma 7.4.\vglue4pt
(iii) We used a rather crude method, namely the identity (6.6) to estimate
the Maclaurin series coefficients of $G_{q}$. It is possible to obtain a far
finer estimation of the coefficients of $G_{q}$. Indeed in \cite
{Lubinsky1999-2}, it was shown that for $q$ such as in (7.2), we have 
\begin{equation}
\left| \log \left| \left( q;q\right) _{n}\right| \right| =O\left( \log
n\right) ,\qquad n\rightarrow \infty ,
\end{equation}
instead of the geometric estimate (7.5). However the problem is that despite
the slower growth in $n$, the size of the constant in the order term in
(7.18) is so large as to render it useless except for large $n$.

\section{Conclusions}

While this paper answers one form of the Baker-Gammel-Wills Conjecture, by
constructing an example in which all subsequences of $\left\{ \left[
n/n\right] \right\} _{n=1}^{\infty }$ have limit points of poles where the
underlying function $H_{q}$ does not, it suggests many more questions. Some
of these are specific to $H_{q}$, but most have a more general flavour.

No matter which $q$ we choose, Worpitzky's theorem ensures that even\break the full
sequence $\left\{ \mu _{n}/\nu _{n}\right\} _{n=1}^{\infty }$ converges
uniformly in compact subsets of\break $\left\{ z:\left| z\right| <\frac{1}{4}%
\right\} $. Can we always obtain uniform convergence of a subsequence near 0?%
\newline

\demo{Problem {\rm 8.1}} {\it Let $f$  be analytic in the unit ball.}
\vglue4pt
(I) {\it Does there exist a neighbourhood of  $0$  and a
subsequence of $\left\{ \left[ n/n\right] \right\} _{n=1}^{\infty }$  converging uniformly to $f$  some ball
centre  $0${\rm ?} If there is such a ball{\rm ,} can its radius be made independent of}  $f$?
\vglue4pt

(II) {\it Does there exist for each point  $a$ in the unit ball{\rm ,}
a neighbourhood of~$a${\rm ,} and a corresponding subsequence of  $%
\left\{ \left[ n/n\right] \right\} _{n=1}^{\infty }$  converging
uniformly to $f$  in that neighbourhood}\/?\vglue4pt

We note that if there is such a ball centre $0$ as in (I), with radius
independent of $f$, then it would imply that for entire functions $f$, there
is a subsequence converging uniformly in all compact subsets of the plane.
This would follow by consideration of $f\left( z/\varepsilon \right) $ with $%
\varepsilon $ small enough, rather than $f\left( z\right) $, and use of
invariance properties for Pad\'{e} approximants. This idea was  used
by\break Buslaev, Gon\v car and Suetin in \cite{Buslaevetal1984}, in their resolution
of the Baker-Graves-Morris conjecture for columns of the Pad\'{e} table.

Another form of the Problem 8.1 involves a weaker form of convergence,
namely convergence in capacity. For the definition and properties of
logarithmic capacity (cap), see for example \cite{Landkof1972}, \cite
{SaffTotik1997}. It is known that for functions analytic in the unit ball,
the full diagonal sequence need not converge in capacity, even in any
neighbourhood of $0$ \cite{Lubinsky1980}, \cite{Lubinsky1983}, \cite
{Rakhmanov1981}. Can a subsequence converge in capacity?\enddemo

{\it Problem} 8.2. {\it Let  $f$  be meromorphic in the unit ball and analytic at}  $%
0 $.\vglue4pt 
(I)  {\it Does there exist a subsequence  $\left\{ \left[ n/n\right]
\right\} _{n\in {\cal S}}$  that converges in capacity to $f${\rm ?} More precisely{\rm ,} given  $0<r<1$  and  $\varepsilon
>0${\rm ,} we want 
$$
{\rm cap}\{z:\left| z\right| \leq r\hbox{ and }\left| f-\left[ n/n\right] \right|
\left( z\right) >\varepsilon \}\rightarrow 0,n\rightarrow \infty ,n\in 
{\cal S}.
$$
If this is not possible for each  $r${\rm ,} is it possible for
some} $r>0?$ 
\vglue4pt
(II) {\it Does there exist a subsequence $\left\{ \left[ n/n\right]
\right\} _{n\in {\cal S}}$ that converges at a given point to  $%
f ${\rm ?} More precisely{\rm ,} we want 
$$
\liminf_{n\rightarrow \infty }\left| \left[ n/n\right] \left( z\right)
-f\left( z\right) \right| =0,
$$
in the unit ball{\rm ,} except at poles of}  $f$.\vglue4pt

We note that problem (I) was raised by H. Stahl in \cite{Stahl1997-2}. The author
recalls that it was also raised by the Russian school of A.\ A.\  
Gon\v car, but
he cannot trace the reference.

Of course, problems about subsequences are difficult to resolve, for if
there is failure, one has to prove this failure for every subsequence. This
is even more difficult when a weak convergence concept, such as convergence
in capacity, is involved.

Is there anything positive that can be said about full diagonal sequences $%
\left\{ \left[ n/n\right] \right\} _{n=1}^{\infty }$ for functions analytic
in the unit ball? Can we find a property weaker than convergence in
capacity, satisfied by the full diagonal sequence? A start in this direction
was made in \cite{Lubinsky2001}, where it is shown that for large $n$, $%
\left[ n/n\right] $ provides good uniform approximation on at least $1/8$ of
the circles centre $0$ within the unit ball.
\enddemo
We now turn to questions concerning the Rogers-Ramanujan function 
$$
G_{q}\left( z\right) =\sum_{j=0}^{\infty }\frac{q^{j^{2}}}{\left( q;q\right)
_{j}}z^{j}.
$$
Recall that the $q$-exponential functions are 
\begin{eqnarray*}
e_{q}\left( z\right) &=&\sum_{j=0}^{\infty }\frac{1}{\left( q;q\right) _{j}}%
z^{j}; \\
E_{q}\left( z\right) &=&\sum_{j=0}^{\infty }\frac{q^{j(j-1)/2}}{\left(
q;q\right) _{j}}z^{j}.
\end{eqnarray*}
Now in $\left\{ z:\left| z\right| <R\left( q\right) \right\} $, a direct
calculation shows that 
$$
e_{q}\left( z\right) E_{q}\left( -z\right) =1
$$
and hence $e_{q}$ and $E_{q}$ have no zeros in that ball. In contrast,
Theorem 2.3 shows that as $q$ approaches a primitive $\ell ^{\rm th}$ root of
unity, the total multiplicity of zeros of $G_{q}$ in a neighbourhood of $%
\left( -1/4\right) ^{-1/\ell }$ approaches $\infty $.

This suggests studying the structure of zeros of $G_{q}$:
\demo{Problem {\rm  8.3}} (i) {\it Investigate the structure of zeros of $G_{q}$ when  $%
\left| q\right| =1$  and}  $R\left( q\right) >0$.
\vglue4pt
(ii) {\it Moreover{\rm ,} investigate whether {\rm -} as seems likely {\rm -} the asymmetry
property {\rm (7.1)} holds for every such}  $q$.\vglue4pt
(iii) {\it Obtain a proof of the asymmetry property {\rm (7.1)} without the use of a
numerical package}.\enddemo

One obvious question is whether our counterexample $H_{q}$ can be used to
provide a counterexample to the Baker-Gammel-Wills Conjecture restricted to
functions $f$ analytic in the unit ball. We can show that for $q$ of (1.8),
(1.9), $G_{q}\left( qz\right) $ has three simple zeros of smallest modulus,
say, 
$$
0<\left| z_{1}\right| <\left| z_{2}\right| <\left| z_{3}\right| <1
$$
and these are the poles of $H_{q}$ of smallest modulus. Let 
\begin{eqnarray*}
g\left( z\right) &:=&H_{q}\left( z\right) (z-z_{1})\left( z-z_{2}\right) ; \\
f\left( z\right) &:=&g\left( zz_{3}\right) .
\end{eqnarray*}
Then $f$ is analytic in the unit ball. 
\demo{Problem {\rm 8.4}} 
 {\it Show that no subsequence of $\left\{ \left[ n/n\right] \right\}
_{n=1}^{\infty }$ to  $f$  can converge uniformly in all
compact subsets of the unit ball.}\enddemo

{\it Acknowledgements}.
The author began the study of $H_{q}$ in 1995, but faltered on the proof of
the asymmetry property (7.1). Arnold Knopfmacher helped the author with
Mathematica pictures for some values of $q$, recorded in \cite
{KnopfmacherLubinsky2001}. These convinced the author that the asymmetry
property (7.1) must hold, although the final choice of $q$ for which a proof
was possible turned out to be different from those in \cite
{KnopfmacherLubinsky2001}. The author thanks Arnold Knopfmacher for those
pictures.

A second acknowledgement must go to the research funding that has supported
the search for a counterexample, which the author began in his doctoral
studies in 1977. At the time, he promised AECI Industries that their
doctoral scholarship would be used to find a counterexample. That promise
has been kept, some twenty years after their generous funding ended. Along
the way research support has been provided by a Lady Davis Fellowship at the
Technion in Haifa, and by the following institutions: the National Research
Institute for Mathematical Sciences in Pretoria; the Centre for Constructive
Mathematics of the University of South Florida; Witwatersrand University in
Johannesburg; the National Research Foundation (formerly FRD) in Pretoria,
and most recently, Georgia Institute of Technology.

Finally the author thanks the referee for requesting error estimates that
led to Lemma 7.4, and for some historical insights. 

\demo{Note added in proof}  V. I. Buslaev has obtained an example of a function analytic in the unit ball for which the
Baker-Gammel-Wills conjecture fails.  His function is algebraic and also negates a conjecture of H. Stahl.  See
\begin{itemize}
\item[{[1]}] V. I. Buslaev, Simple counterexample to the Baker-Gammel-Wills conjecture, {\it East Journal on
Approximations} {\bf 4} (2001), 515--517.

\item[{[2]}] V. I. Buslaev, On the Baker-Gammel-Wills conjecture in the theory of Pad\'e approximants, {\it Mathematics
Sbornik}, to appear.
\end{itemize}

\end{document}